\documentclass{gtart_h}

\def\ifplaintex{\expandafter\ifx\csname documentclass\endcsname\relax}

\def\gtp{{\mathsurround=0pt\it $\cal G\mskip-2mu$eometry \&\ 
$\cal T\!\!$opology $\cal P\!$ublications}}  

\def\recd{{\small Received:\qua\receiveddate\ifx\reviseddate\relax
\else\qquad Revised:\qua\reviseddate\fi\par}} 

\def\lognumber#1{\def\thelognumber{#1}}
\def\volumenumber#1{\def\thevolumenumber{#1}}
\def\volumeyear#1{\def\thevolumeyear{#1}}
\def\papernumber#1{\def\thepapernumber{#1}}
\def\pagenumbers#1#2{\def\startpage{#1}\def\finishpage{#2}}
\def\published#1{\def\publishdate{#1}}

\def\received#1{\def\receiveddate{#1}}
\def\revised#1{\def\reviseddate{#1}}
\def\accepted#1{\def\accepteddate{#1}}
\def\asciititle#1{\def\theasciititle{#1}}

\def\asciiaddress#1{\def\theasciiaddress{#1}}
\def\asciiemail#1{\def\theasciiemail{#1}}

\long\def\asciiabstract#1{\long\def\theasciiabstract{#1}}

\let\\\par\let\thelognumber\relax\let\thevolumenumber\relax
\let\thepapernumber\relax\let\thevolumeyear\relax\let\startpage\relax
\let\finishpage\relax\let\publishdate\relax\let\receiveddate\relax
\let\reviseddate\relax\let\accepteddate\relax\let\theasciititle\relax
\let\theasciiauthors\relax\let\theasciiaddress\relax
\let\theasciiabstract\relax

\let\theasciiemail\relax

\ifplaintex
\font\logobig=cmssbx10 scaled 3836
\font\logomed=cmssbx10 scaled 2557
\else
\font\logobig=cmssbx10 scaled 4200
\font\logomed=cmssbx10 scaled 2800
\fi

\long\def\makeagttitle{   
\count0=\startpage
\agt\hfill      
\hbox to 45truept{\vbox to 0pt{\vglue -13truept{\logomed A\kern -.37em{\logobig 
T}\kern -.38em G}\vss}\hss}
\break
{\small Volume \thevolumenumber\ (\thevolumeyear)
\startpage--\finishpage\nl
Published: \publishdate}

\vglue .25truein

{\parskip=0pt\leftskip 0pt plus
1fil\def\\{\par\smallskip}{\Large\bf\thetitle}\par\medskip} \vglue
0.05truein

{\parskip=0pt\leftskip 0pt plus 1fil\def\\{\par}{\sc\theauthors}
\par\medskip}%
 
\vglue 0.03truein 

{\small\leftskip 25truept\rightskip 25truept{\bf Abstract}\stdspace\theabstract

{\bf AMS Classification}\stdspace\theprimaryclass
\ifx\thesecondaryclass\relax\else; \thesecondaryclass\fi\par
{\bf Keywords}\stdspace \thekeywords\par}\vglue 7truept

}   

\ifplaintex
\font\phead=cmsl9 scaled 950
\font\pnum=cmbx10 scaled 913
\font\pfoot=cmsl9 scaled 950
\headline{\vbox to 0pt{\vskip -4.5mm\line{\small\phead\ifnum
\count0=\startpage ISSN 1472-2739 (on-line) 1472-2747 (printed)
\hfill {\pnum\folio}\else\ifodd\count0\def\\{ }%
\ifx\theshorttitle\relax\thetitle\else\theshorttitle\fi\hfill{\pnum\folio}
\else\def\\{ and }{\pnum\folio}\hfill\ifx\theshortauthors\relax\theauthors
\else\theshortauthors\fi\fi\fi}\vss}}
\footline{\vbox to 0pt{\vglue 0mm\line{\small\pfoot\ifnum\count0=\startpage
\copyright\ \gtp\hfill\else
\agt, Volume \thevolumenumber\ (\thevolumeyear)\hfill\fi}\vss}}
\else
\font=cmsl9 scaled 1050
\font\lnum=cmbx10 
\font=cmsl9 scaled 1050
\makeatletter
\def\@oddhead{{\small\lhead\ifnum\count0=\startpage ISSN 1472-2739 
(on-line) 1472-2747 (printed)\hfill {\lnum\number\count0}\else\ifodd\count0
\def\\{ }\ifx\theshorttitle\relax \thetitle \else\theshorttitle\fi\hfill
{\lnum\number\count0}\else\def\\{ and }{\lnum\number\count0}
\hfill\ifx\theshortauthors\relax 
\theauthors\else\theshortauthors\fi\fi\fi}}\def\@evenhead{\@oddhead}
\def\@oddfoot{\small\lfoot\ifnum\count0=\startpage\copyright\ \gtp\hfill\else
\agt, Volume \thevolumenumber\ (\thevolumeyear)\hfill\fi}
\def\@evenfoot{\@oddfoot}
\makeatother
\fi
\let\maketitlepage\makeagttitle
\let\makeshorttitle\maketitlepage
\let\maketitle\maketitlepage

\newwrite\gtoutfile
\long\gdef\makeheadfile{  
{\def\\{, }\def\s{ }
\immediate\openout\gtoutfile head.xxx
\immediate\write\gtoutfile{Proxy-for: \ifx\theasciiauthors\relax
\theauthors\else\theasciiauthors\fi\s<\ifx\theasciiemail\relax\theemail\else\theasciiemail\fi>}
\immediate\write\gtoutfile{\noexpand\\}
\immediate\write\gtoutfile{Authors: \ifx\theasciiauthors\relax
\theauthors\else\theasciiauthors\fi}
{\def\\{ }\immediate\write\gtoutfile{Title: \ifx\theasciititle\relax
\thetitle\else\theasciititle\fi}}
\immediate\write\gtoutfile{Subj-class: GT or SG, GR etc}
\immediate\write\gtoutfile{MSC-class: \theprimaryclass\ifx\thesecondaryclass\relax\else, \thesecondaryclass\fi}
\immediate\write\gtoutfile{Journal-ref: Algebr. Geom. Topol. \thevolumenumber\s
(\thevolumeyear) \startpage-\finishpage}
\immediate\write\gtoutfile{Comments: Published by Algebraic and
Geometric Topology at}
\immediate\write\gtoutfile{\s\s\s  http://www.maths.warwick.ac.uk/agt/AGTVol\thevolumenumber/agt-\thevolumenumber-\thepapernumber.abs.html}
\immediate\write\gtoutfile{\noexpand\\}
\immediate\write\gtoutfile{}
\ifx\theasciiabstract\relax
\immediate\write\gtoutfile{\theabstract}\else
\immediate\write\gtoutfile{\theasciiabstract}\fi
\immediate\write\gtoutfile{}
\immediate\write\gtoutfile{\noexpand\\}
\immediate\write\gtoutfile{}
\immediate\closeout\gtoutfile}}  

\def\maketitlepage{\makeagttitle\makeheadfile}
\let\makeshorttitle\maketitlepage
\let\maketitle\maketitlepage

\lognumber{21}
\volumenumber{5}
\volumeyear{2005}
\papernumber{21}
\pagenumbers{463}{507}
\received{11 January 2005} 
\revised{13 April 2005}
\accepted{29 April 2005}
\published{30 May 2005}

\usepackage{amssymb,amsmath,graphicx}

\newtheorem{theorem}{Theorem}[section]
\newtheorem{lemma}[theorem]{Lemma}

\newtheorem{corollary}[theorem]{Corollary}
\newtheorem{proposition}[theorem]{Proposition}
\newtheorem{claim}[theorem]{Claim}
\theoremstyle{definition}

\numberwithin{equation}{section}

\begin{document}

\title[Maximal distance between toroidal Dehn fillings]{On hyperbolic 
$3$-manifolds realizing the maximal\\distance between toroidal Dehn fillings}
\asciititle{On hyperbolic 3-manifolds realizing the maximal distance between toroidal Dehn fillings}

\authors{Hiroshi Goda\\Masakazu Teragaito}

\address{Department of Mathematics, Tokyo University of Agriculture 
and Technology\\Koganei, Tokyo 184-8588, Japan}
\secondaddress{Department of Mathematics and Mathematics Education, Hiroshima University\\1-1-1 Kagamiyama, Higashi-hiroshima, Japan 739-8524}
\asciiaddress{Department of Mathematics, Tokyo University of Agriculture 
and Technology\\Koganei, Tokyo 184-8588, Japan\\and\\Department of 
Mathematics and Mathematics Education, Hiroshima University\\1-1-1 
Kagamiyama, Higashi-hiroshima, Japan 739-8524}
\gtemail{\mailto{goda@cc.tuat.ac.jp}{\rm\qua and\qua}\mailto{teragai@hiroshima-u.ac.jp}}
\asciiemail{goda@cc.tuat.ac.jp, teragai@hiroshima-u.ac.jp}

\begin{abstract}
For a hyperbolic $3$-manifold $M$ with a torus boundary component,
all but finitely many Dehn fillings on the torus component yield hyperbolic $3$-manifolds.
In this paper, we will focus on the situation where 
$M$ has two exceptional Dehn fillings, both of which yield toroidal manifolds.
For such situation, Gordon gave an upper bound for the distance between two slopes of Dehn fillings.
In particular, if $M$ is large, then the distance is at most $5$.
We show that this upper bound can be improved by $1$ for a broad class of large manifolds.
\end{abstract}

\asciiabstract{%
For a hyperbolic 3-manifold M with a torus boundary component, all but
finitely many Dehn fillings on the torus component yield hyperbolic
3-manifolds.  In this paper, we will focus on the situation where M
has two exceptional Dehn fillings, both of which yield toroidal
manifolds.  For such situation, Gordon gave an upper bound for the
distance between two slopes of Dehn fillings.  In particular, if M is
large, then the distance is at most 5.  We show that this upper bound
can be improved by 1 for a broad class of large manifolds.}


\maketitle

\section{Introduction}\label{sec:1}

Let $M$ be a hyperbolic $3$-manifold with a torus boundary component $T_0$.
A \textit{slope\/} on $T_0$ is the isotopy class of an essential simple closed curve on $T_0$.
For a slope $\gamma$ on $T_0$, the manifold obtained by \textit{$\gamma$-Dehn filling\/} is $M(\gamma)=M\cup V_\gamma$, where
$V_\gamma$ is a solid torus, glued to $M$ along $T_0$ in such a way that $\gamma$ bounds a meridian disk in $V_\gamma$.
If $M(\gamma)$ is not hyperbolic, then $\gamma$ is called an \textit{exceptional slope}.
By Thurston's hyperbolic Dehn surgery theorem, 
the number of exceptional slopes is finite.
If $M(\gamma)$ fails to be hyperbolic, then it either (1) contains an essential sphere, disk, annulus or torus; or
(2) contains a Heegaard sphere or torus; or (3) is a Seifert fibered manifold over the sphere with three
exceptional fibers; or (4) is a counterexample to the geometrization conjecture
(see \cite{Go4}).

Suppose that there are two slopes $\alpha$ and $\beta$ such that $M(\alpha)$ and $M(\beta)$ are
toroidal, that is, contain essential tori.
The \textit{distance\/} $\Delta(\alpha,\beta)$ between them is their minimal geometric intersection number.
Then Gordon \cite{Go2} shows $\Delta=\Delta(\alpha,\beta)\le 8$, and
there are only four manifolds $W(-1), W(5), W(5/2), W(-2)$ with $\Delta\ge 6$.
Here, $W(p/q)$ is obtained by $p/q$-filling on one boundary torus of the Whitehead link exterior $W$ in the usual way.
In particular, these manifolds are each $\mathbb{Q}$-homology $S^1\times D^2$, and the boundary is a single torus.
Following Wu \cite{W2},
let us say that $M$ is \textit{large\/} if $H_2(M,\partial M-T_0)\ne 0$.
Note that $M$ is not large if and only if $M$ is a $\mathbb{Q}$-homology $S^1\times D^2$ or a $\mathbb{Q}$-homology $T^2\times I$. 
Hence, $M$ is large if $\partial M$ is not a union of at most two tori.
In \cite[Question 4.2]{Go4}, Gordon asks if there is a large hyperbolic manifold with toroidal fillings at distance $5$.
In this direction, \cite[Theorem 3.1]{BGZ} shows that if $\partial M$ is a single torus and the first betti number $\beta_1(M)\ge 3$ then
the distance between two toroidal fillings is at most $4$.
As stated in \cite[Remark 3.15]{BGZ}, their argument also works for $M$ whose boundary consists of at least $4$ tori.

The purpose of this paper is to show that a broad class of large manifolds cannot admit two toroidal fillings at distance $5$. 

\begin{theorem}\label{thm:main}
Let $M$ be a hyperbolic $3$-manifold with a torus boundary component $T_0$ and 
suppose that there are two slopes $\alpha$, $\beta$ on $T_0$ such that $M(\alpha)$ and $M(\beta)$ are toroidal.
If $\Delta(\alpha,\beta)=5$, then $\partial M$ consists of at most two tori.
\end{theorem}

This is sharp in the sense that there are hyperbolic $3$-manifolds whose boundary is a single or two tori with
two toroidal fillings at distance $5$.
For example, the exterior of the $(-2,3,7)$-pretzel knot in $S^3$ is hyperbolic
and there are two toroidal slope $16$ and $37/2$.
The Whitehead sister link ($(-2,3,8)$-pretzel link) exterior gives such an example with two torus boundary components.
Also, Theorem \ref{thm:main} can be regarded as the first step to determine which hyperbolic $3$-manifolds admit
two toroidal slopes of distance $5$.
Part of the proof of Theorem \ref{thm:main} consists of carrying over the argument of \cite{T},
where we treated the case where $M$ is the exterior of a hyperbolic knot in $S^3$, to the present context.
Hence we assume the familiarity with \cite{T}.

Theorem \ref{thm:main} gives also a partial answer to \cite[Question 5.2]{Go4} which asks if there is a hyperbolic manifold
whose boundary consists of three tori, having two toroidal fillings at distance $4$ or $5$.
Combining with known facts \cite{Go4}, we have the following.

\begin{corollary}
If $M$ is a hyperbolic $3$-manifold whose boundary is a union of more than two tori, then 
for any fixed boundary torus component $T_0$ of $M$, any two exceptional slopes of $M$ on $T_0$
have distance at most $4$.
\end{corollary}

To prove Theorem \ref{thm:main}, we need to consider the situation where
either $M(\alpha)$ or $M(\beta)$ contains a Klein bottle.
Such a phenomenon often happens in the literature \cite{GL2,GL3,LOT}.

\begin{theorem}\label{thm:klein}
Let $M$ be a hyperbolic $3$-manifold with a torus boundary component $T_0$ and 
suppose that there are two slopes $\alpha$, $\beta$ on $T_0$ such that $M(\alpha)$ contains a Klein bottle and $M(\beta)$ is toroidal.
If $\Delta(\alpha,\beta)\ge 5$, then $\partial M$ consists of at most two tori.
\end{theorem}

In Section \ref{sec:pre}, we prepare some general lemmas about a pair of graphs coming from
intersections of two essential tori.
Sections \ref{sec:generic}--\ref{sec:st2} treat the case where two toroidal manifolds contain no Klein bottle.
Finally, we consider the case where either contains a Klein bottle in Section \ref{sec:klein}--\ref{sec:generic-kb}.
Section \ref{sec:reduced} contains the results about a reduced graph on a Klein bottle, which we need for Section \ref{sec:generic-kb}.

\section{Preliminaries}\label{sec:pre}

Let $M$ be a hyperbolic $3$-manifold with a torus boundary component $T_0$ and 
suppose that there are two slopes $\alpha$, $\beta$ on $T_0$ such that $M(\alpha)$ and $M(\beta)$ are toroidal.
We assume that $\Delta=\Delta(\alpha,\beta)=5$ until the end of Section \ref{sec:st2}.
Then $M(\alpha)$ and $M(\beta)$ are irreducible by \cite{O,W}.

Let $\widehat{S}$ be an essential torus in $M(\alpha)$.
We may assume that $\widehat{S}$ meets the attached solid torus $V_\alpha$ in $s$ meridian disks
$u_1,u_2,\dots,u_s$, numbered successively along $V_\alpha$, and
that $s$ is minimal over all choices of $\widehat{S}$.
Let $S=\widehat{S}\cap M$.
Then $S$ is a punctured torus properly embedded in $M$ with $s$ boundary components $\partial_iS=\partial u_i$,
each of which has slope $\alpha$.
By the minimality of $s$, $S$ is incompressible and boundary-incompressible in $M$.
Similarly, we choose an essential torus $\widehat{T}$ in $M(\beta)$ which intersects the
attached solid torus $V_\beta$ in $t$ meridian disks $v_1,v_2,\dots,v_t$, numbered successively along $V_\beta$, where $t$ is minimal as above.
Then we have another incompressible and boundary-incompressible punctured torus $T=\widehat{T}\cap M$,
which has $t$ boundary components $\partial_jT=\partial v_j$.
Notice that $s$ and $t$ are non-zero.

We may assume that $S$ intersects $T$ transversely.
Then $S\cap T$ consists of arcs and circles.
Since both surfaces are incompressible, we can assume that no circle component of $S\cap T$
bounds a disk in $S$ or $T$.
Furthermore, it can be assumed that $\partial_iS$ meets $\partial_jT$ in $5$ points for any pair of $i$ and $j$.

As seen in \cite{GW},
we can choose a meridian-longitude pair $m, l$ on $T_0$ so that $\alpha=m$, and $\beta=dm+5l$ for some $d=1,2$.
This number $d$ is called the \textit{jumping number\/} of $\alpha$ and $\beta$.

\begin{lemma}\label{lem:jumping}
Let $a_1,a_2,a_3,a_4,a_5$ be the points of $\partial_iS\cap\partial_jT$, numbered so that they appear successively on $\partial_iS$.
If $d$ is the jumping number of $\alpha$ and $\beta$,
then these points appear in the order of $a_d,a_{2d},a_{3d},a_{4d},a_{5d}$ on $\partial_jT$ in some direction.
In particular, 
if $d=1$, then
two points of $\partial_iS\cap\partial_jT$ are successive on $\partial_iS$ if and only if they are successive in $\partial_jT$,
and if $d=2$, then
two points of $\partial_iS\cap\partial_jT$ are successive on $\partial_iS$ if and only if they are not successive in $\partial_jT$.
\end{lemma}

\begin{proof}
See \cite[Lemma 2.10]{GW}.
\end{proof}

Let $G_S$ be the graph on $\widehat{S}$ consisting of the $u_i$ as (fat) vertices, and the arc components
of $S\cap T$ as edges.
Each vertex of $G_S$ is given a sign according to whether the core of $V_\alpha$ passes $\widehat{S}$ from the
positive side or negative side at this vertex.
Define $G_T$ on $\widehat{T}$ similarly.
Throughout the paper, two graphs on a surface are considered to be equivalent if there is a homeomorphism of the surface
carrying one graph to the other.
Note that $G_S$ and $G_T$ have no trivial loops, since $S$ and $T$ are boundary-incompressible.

For an edge $e$ of $G_S$ incident to $u_i$,
the endpoint of $e$ is labelled $j$ if it is in $\partial u_i\cap \partial v_j=\partial_iS\cap\partial_jT$.
Similarly, label the endpoints of each edge of $G_T$.
Thus the labels $1,2,\dots,t$ (resp.\ $1,2,\dots,s$) appear in order around each vertex of $G_S$ (resp.\ $G_T$) repeated $5$ times.
Each vertex $u_i$ of $G_S$ has degree $5t$, and each $v_j$ of $G_T$ has degree $5s$.

Let $G=G_S$ or $G_T$.
An edge of $G$ is a \textit{positive\/} edge if it connects vertices of the same sign.
Otherwise it is a \textit{negative\/} edge.
Possibly, a positive edge is a loop.
An endpoint of a positive (resp.\ negative) edge around a vertex is called a \textit{positive} (resp.\ \textit{negative\/}) \textit{edge endpoint\/}.
We denote by $G^+$ the subgraph of $G$ consisting of all vertices and positive edges of $G$.

If an edge $e$ of $G_S$ is incident to $u_i$ with label $j$, then it is called a \textit{$j$-edge at} $u_i$.
Then $e$ is also an $i$-edge at $v_j$ in $G_T$.
If $e$ has labels $j_1,j_2$ at its endpoints, then $e$ is called a \textit{$\{j_1,j_2\}$-edge}.
An $\{i,i\}$-edge is said to be \textit{level\/}.

A cycle in $G$ consisting of positive edges is a \textit{Scharlemann cycle\/} if it bounds a disk face of $G$ and all edges
in the cycle are $\{i,i+1\}$-edges for some label $i$.
The number of edges in a Scharlemann cycle is called the \textit{length\/} of the Scharlemann cycle, and
the set $\{i,i+1\}$ is called its \textit{label pair}.
A Scharlemann cycle of length two is called an \textit{$S$-cycle} for short.
For a label $x$, let $G_x$ be the subgraph of $G$ consisting of all vertices and all positive $x$-edges.
Then a disk face of $G_x$ is called an \textit{$x$-face}. 

\begin{lemma}\label{lem:common}
\begin{itemize}
\item[\rm(1)] \textup{(The parity rule)}\qua An edge $e$ is positive in a graph if and only if it is negative in the other graph.
\item[\rm(2)] There is no pair of edges which are parallel in both graphs.
\item[\rm(3)] If $G_S$ \textup{(}resp. $G_T$\textup{)} has a Scharlemann cycle, then $\widehat{T}$ \textup{(}resp. $\widehat{S}$\textup{)} is separating.  
\end{itemize}
\end{lemma}

\begin{proof}
(1)\qua This can be found in \cite{CGLS}.
(2) is \cite[Lemma 2.1]{Go2}. See \cite{BZ0} for (3).
\end{proof}

\begin{proposition}
Either $\widehat{S}$ or $\widehat{T}$ is separating.
\end{proposition}

\begin{proof}
If $G_T$ has more than $t$ positive $x$-edges for some label $x$, then
$G_T$ has an $x$-face, which contains a Scharlemann cycle by \cite{HM}.
Then $\widehat{S}$ is separating by Lemma \ref{lem:common}(3).

Hence we assume that $G_T$ has at most $t$ positive $x$-edges for any label $x$.
This means that any vertex of $G_S$ is incident to at most $t$ negative edges by the parity rule.
Thus any vertex of $G_S$ has at least $4t$ positive edge endpoints, and then 
$G_S^+$ has at least $2st$ edges.
But this implies that $G_S$ has more than $s$ positive $i$-edges for some label $i$.
Then $G_S$ has an $i$-face, containing a Scharlemann cycle.
So $\widehat{T}$ is separating by Lemma \ref{lem:common}(3) again.
\end{proof}

Thus we can assume that $\widehat{S}$ is separating until the end of Section \ref{sec:st2}.
Then $s$ is even.
Let $M(\alpha)=\mathcal{B}\cup_{\widehat{S}}\mathcal{W}$.
Here $\mathcal{B}$ is called the black side of $\widehat{S}$, and $\mathcal{W}$ is the white side.
A Scharlemann cycle is said to be \textit{black\/} (resp. \textit{white\/}) if its face lies in $\mathcal{B}$ (resp. $\mathcal{W}$).

\begin{lemma}\label{lem:parallel-max2}
$G_S$ satisfies the following\textup{:}
\begin{itemize}
\item[\rm(1)] If $\widehat{T}$ is non-separating, then any family of parallel positive edges in $G_S$ contains at most $t/2$ edges.
If $\widehat{T}$ is separating and $t\ge 4$, then any family of parallel positive edges in $G_S$ contains at most $t/2+2$ edges, and 
moreover, if the family contains $t/2+2$ edges, then $t\equiv 0\pmod{4}$, and
$M(\beta)$ contains a Klein bottle. 
\item[\rm(2)] Either any family of parallel negative edges in $G_S$ contains at most $t$ edges, or
all vertices of $G_T$ have the same sign.
\end{itemize}
\end{lemma}

\begin{proof}
(1)\qua If $\widehat{T}$ is non-separating, then $G_S$ cannot contain a Scharlemann cycle by Lemma \ref{lem:common}(3).
Thus any family of parallel positive edges in $G_S$ contains at most $t/2$ edges by \cite[Lemma 2.6.6]{CGLS}.
Assume that $\widehat{T}$ is separating and $t\ge 4$.
By \cite[Lemma 1.4]{W}, any family of parallel positive edges contains at most $t/2+2$ edges.
If the family contains $t/2+2$ edges, then $t\equiv 0\pmod{4}$ by \cite[Corollary 1.8]{W}. 
In this case, the family contains two $S$-cycles $\rho_1$ and $\rho_2$ with disjoint label pairs.
Let $\{k_i,k_i+1\}$ be the label pair of $\rho_i$ and let $D_i$ be the disk face bounded by $\rho_i$ for $i=1,2$.
Let $H_i$ be the part of $V_\beta$ between $v_{k_i}$ and $v_{k_i+1}$.
Then shrinking $H_i$ into its core in $H_i\cup D_i$ gives a M\"{o}bius band $B_i$ whose boundary is the loop on $\widehat{T}$
formed by the edges of $\rho_i$. 
In particular, $\partial B_i$ is essential on $\widehat{T}$ \cite[Lemma 3.1]{GL2}.
Hence $\partial B_1$ and $\partial B_2$ are disjoint, and so they bound an annulus $A$ on $\widehat{T}$.
Then the union $B_1\cup A\cup B_2$ is a Klein bottle in $M(\beta)$.

(2)\qua If $t=1$, then the second conclusion holds.
If $t=2$, then $G_T$ has only two parallelism classes of loops \cite[Lemma 5.2]{Go2}.
Hence if two vertices of $G_T$ have opposite signs, then at most two negative edges can be parallel in $G_S$ by Lemma \ref{lem:common}(2).
See \cite[Lemma 2.3(1)]{GW} for $t>2$.
\end{proof}

\begin{lemma}\label{lem:parallel-max}
$G_T$ satisfies the following\textup{:}
\begin{itemize}
\item[\rm(1)] If $s\ge 4$, then any family of parallel positive edges in $G_T$ contains at most $s/2+2$ edges.
Moreover, if the family contains $s/2+2$ edges, then $s\equiv 0\pmod{4}$, and
$M(\alpha)$ contains a Klein bottle. 
\item[\rm(2)] Any family of parallel negative edges in $G_T$ contains at most $s$ edges.
\end{itemize}
\end{lemma}

\begin{proof}
This can be proved by the same argument as in the proof of Lemma \ref{lem:parallel-max2}.
\end{proof}

For a graph $G$ on a surface,
$\overline{G}$ denotes the \textit{reduced graph} of $G$ obtained by amalgamating each family of parallel edges into a single edge.
For an edge $e$ of $\overline{G}$, the \textit{weight\/} of $e$ is the number of edges in the corresponding family of parallel edges in $G$.

\section{Generic case}\label{sec:generic}

The proof of Theorem \ref{thm:main} occupies Sections \ref{sec:generic}--\ref{sec:st2}.
The case where either $M(\alpha)$ or $M(\beta)$ contains a Klein bottle will be treated from Section \ref{sec:klein}.
Hence we assume that neither $M(\alpha)$ nor $M(\beta)$ contains a Klein bottle in the following $5$ sections.
This section treats the case where $s\ge 4$ and $t\ge 3$.

\begin{lemma}\label{lem:noklein}
\begin{itemize}
\item[\rm(1)] Any family of mutually parallel positive edges in $G_S$ \textup{(}resp.\ $G_T$\textup{)} contains
at most $t/2+1$ \textup{(}resp.\ $s/2+1$\textup{)} edges.
\item[\rm(2)] Neither $G_S$ nor $G_T$ contains two $S$-cycles with disjoint label pairs.
\end{itemize}
\end{lemma}

\begin{proof}
(1) follows from Lemmas \ref{lem:parallel-max2}(1) and \ref{lem:parallel-max}(1).

(2)\qua If $G_S$, say, contains two $S$-cycles with disjoint label pairs, then
$M(\beta)$ contains a Klein bottle as in the proof of Lemma \ref{lem:parallel-max2}.
\end{proof}


Under the existence of Lemma \ref{lem:noklein},
we can carry over the arguments from Lemma 4.1 to 4.13 of \cite{T}.
(In the proof of Lemma 4.12 of \cite{T}, we need to add the case where $t=3$, but it is obvious.)
Hence we have $s=4$ or $6$.
To eliminate these remaining cases, we have to modify the arguments in \cite{T},
because $\widehat{T}$ is possibly non-separating, and the jumping number is one or two in the present context.
(In \cite{T}, both tori were separating and the jumping number between the slopes was one.)

\begin{proposition}
$s=6$ is impossible.
\end{proposition}

\begin{proof}
By \cite[Lemma 4.13]{T}, $\overline{G}_S^+$ consists of two components, each of which has three vertices.
Also, $\overline{G}_S^+$ has a good vertex $u_i$ of degree three, and $t\le 6$ (see the first paragraph of the proof of \cite[Proposition 4.14]{T}).

Assume $t=6$.
If $u_i$ has more than $18(=3t)$ negative edge endpoints in $G_S$, then some label appears four times there.
This implies $s=4$ by \cite[Lemma 4.7]{T}.
Hence it suffices to consider the case where $u_i$ is incident to three families of $4(=s/2+1)$ parallel positive edges.
Then there are just $18$ negative edge endpoints successively at $u_i$.
Thus any label $j$ appears three times among there.
In $G_T$, the vertex $v_j$ is incident to three positive $i$-edges.
No two of them are parallel by Lemma \ref{lem:parallel-max}(1), and so
$v_j$ is incident to three families of $4$ parallel positive edges and three families of $6$ parallel negative edges.
Notice that each of the three families of positive edges contains an $i$-edge with label $i$ at $v_j$.
But it is easy to see that such labeling is impossible around $v_j$.
The case $t=4$ is similar to this case.

If $t=5$, then any positive edge at $u_i$ has weight at most two, since $G_S$ cannot contain a Scharlemann cycle.
Thus $u_i$ is incident to at least $19$ negative edges successively in $G_S$.
Then some label appears $4$ times among negative edge endpoints.
This implies $s=4$ by \cite[Lemma 4.7]{T}.
The case $t=3$ is similar to this.
\end{proof}

\begin{figure}[htb]
\begin{center}
\includegraphics*[scale=0.26]{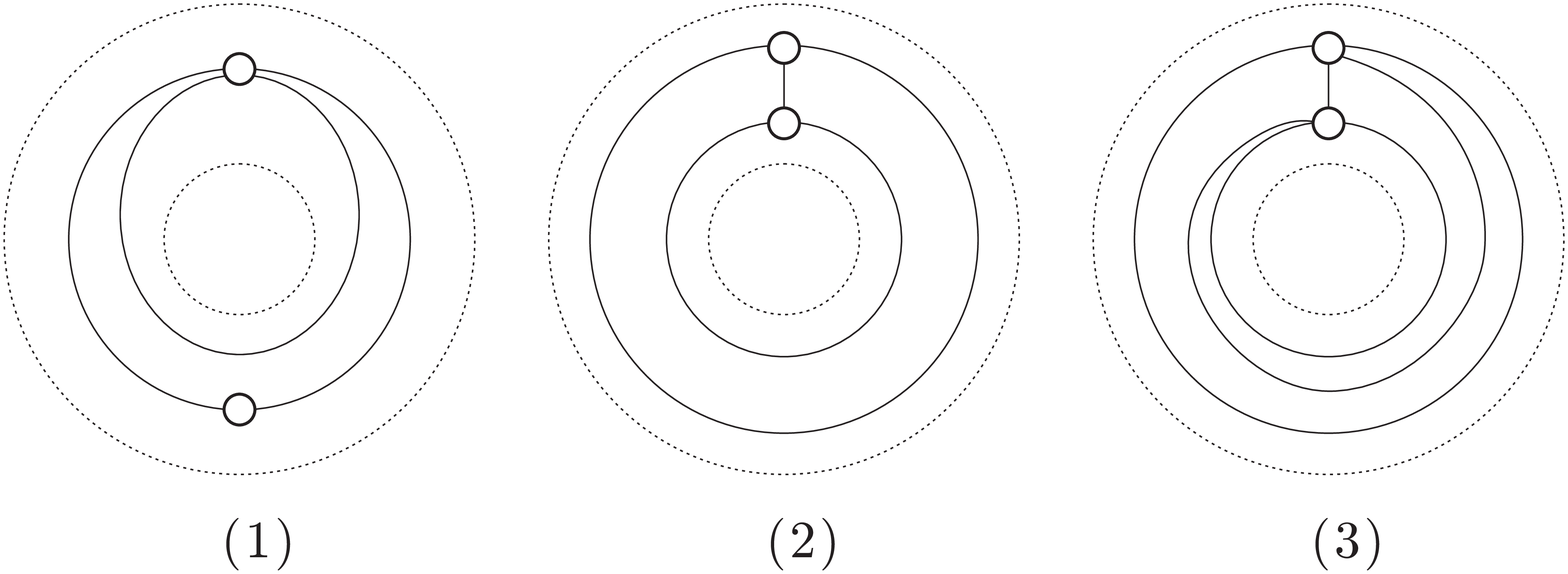}
\caption{}\label{fig:s4}
\end{center}
\end{figure}

Therefore we have $s=4$.
Then $\overline{G}_S^+$ consists of two components, each of which has the form of Figure \ref{fig:s4}(1), (2) or (3) by \cite[Lemmas 4.8, 4.11]{T}.

\nocolon

\begin{lemma}\label{lem:s4-1}
$\overline{G}_S^+$ does not have a component of the form as in Figure \textup{\ref{fig:s4}(1)}.
\end{lemma}

\begin{proof}
Let $\Gamma$ be a component of $\overline{G}_S^+$ as in Figure \ref{fig:s4}(1),
and let $u_i$ be the good vertex of degree two in $\Gamma$.
Since $u_i$ is incident to at most $2(t/2+1)=t+2$ positive edges in $G_S$,
there are at least $4t-2$ negative edge endpoints.
Then some label $j$ appears $4$ times there, because $4t-2>3t$.
In $G_T$, $v_j$ is incident to $4$ positive $i$-edges.
Since no two of them are parallel, $v_j$ is incident to $4$ families of $3$ parallel positive edges and
two families of $4$ parallel negative edges.
Notice that each family of positive edges contains an $i$-edge with label $i$ at $v_j$.
But this is clearly impossible, because both families of negative edges contain $i$-edges.
(Recall that any label appears just $5$ times around a vertex.)
\end{proof}

\begin{lemma}\label{lem:s4-2}
$\overline{G}_S^+$ does not have a component of the form as in Figure \textup{\ref{fig:s4}(2)}.
\end{lemma}

\begin{proof}
Let $\Gamma$ be such a component with a good vertex $u_i$ of degree three.
Assume $t>6$.
Then $u_i$ has at least $5t-3(t/2+1)=7t/2-3>3t$ negative edge endpoints.
Hence some label $j$ appears $4$ times there.
Then the same argument as in the proof of Lemma \ref{lem:s4-1} works.
If $t=3$ or $5$, then $u_i$ has more than $3t$ negative edge endpoints, because $G_S$ cannot contain
a Scharlemann cycle.
Then some label appears $4$ times again, and so it leads to a contradiction.

Finally, the argument in the proof of \cite[Lemma 4.16]{T} works when $t=4$ and $6$.
(Use Lemma \ref{lem:noklein}(2) instead of \cite[Lemma 2.7(2)]{T}.)
\end{proof}

\begin{proposition}
$s=4$ is impossible.
\end{proposition}

\begin{proof}
By Lemmas \ref{lem:s4-1} and \ref{lem:s4-2}, $\overline{G}_S^+$ consists of two components of the form as in Figure \ref{fig:s4}(3).
Notice that any vertex of $\overline{G}_S$ is incident to at most two negative edges.
Let $u$ be a vertex of $G_S$.

First, suppose that $\widehat{T}$ is separating.  Then $t\ge 4$.
Hence any family of parallel negative edges in $G_S$ has at most $t$ edges by Lemma \ref{lem:parallel-max2}.
Thus $u$ has at most $2t$ negative edge endpoints, and then it has at least $3t$ positive edge endpoints.
From $4(t/2+1)\ge 3t$, we have $t=4$.
Then $u$ is incident to three loops and two families of $3$ parallel positive edges, and
so there are two $S$-cycles with disjoint label pairs, which is impossible by Lemma \ref{lem:noklein}.

Hence $\widehat{T}$ is non-separating.
Then $u$ has at most $4\cdot t/2=2t$ positive edge endpoints by Lemma \ref{lem:parallel-max2}.
Hence there are at least $3t$ negative edge endpoints consecutively.
If there are more than $3t$, then some label appears $4$ times there, which leads to a contradiction as in the proof of Lemma \ref{lem:s4-1}.
Thus $u$ has exactly $3t$ negative edge endpoints, and is incident to $4$ families of $t/2$ parallel positive edges.

Let $u'$ be the other vertex of the same component as $u$.
Since $G_S$ cannot contain a Scharlemann cycle, the labeling around $u'$ is uniquely determined by the labeling around $u$.
But then it is clear to see that there is a Scharlemann cycle of length three.
\end{proof}

\section{The case $t=1$}\label{sec:t1}

The reduced graph $\overline{G}_T$ consists of at most three edges by \cite[Lemma 5.1]{Go2}.
We denote the weights of the edges by $(w_1,w_2,w_3)$, and say $G_T\cong G(w_1,w_2,w_3)$ as in \cite{Go2}.
Notice that $G(w_1,w_2,w_3)$ is invariant under any permutations of the $w$'s.

\begin{lemma}
$s=2$.
\end{lemma}

\begin{proof}
If $s\ge 4$, then the vertex of $G_T$ is incident to at most $6(s/2+1)=3s+6$ edges.
From $3s+6\ge 5s$, we have $s\le 3$, a contradiction.
\end{proof}

Thus $G_T$ has exactly five $\{1,2\}$-edges, which are divided into at most three families of mutually
parallel edges.
Since any edge of $G_T$ is positive, all edges of $G_S$ are negative by the parity rule, and they are divided into
at most $4$ classes (see \cite{GL3}).

\begin{lemma}
$G_T\cong G(3,1,1)$.
\end{lemma}

\begin{proof}
If two parallel edges of $G_T$ have the same edge class label, then
these edges are parallel in both $G_S$ and $G_T$.
This is impossible by Lemma \ref{lem:common}.
Hence at most four edges can be parallel in $G_T$.
Then $G_T\cong G(4,1,0)$, $G(3,2,0)$, $G(3,1,1)$ or $G(2,2,1)$.
However all but $G(3,1,1)$ are impossible, because each edge must have labels $1$ and $2$ at its endpoints.
\end{proof}

\begin{proposition}\label{prop:t1case}
$\partial M$ consists of a single torus. 
\end{proposition}

\begin{proof}
If the jumping number is one,
then $G_S$ and $G_T$ are determined as shown in Figure \ref{fig:t1s2-jump1}, where
the correspondence of edges is indicated.

\nocolon\begin{figure}[htb]
\begin{center}
\includegraphics*[scale=0.6]{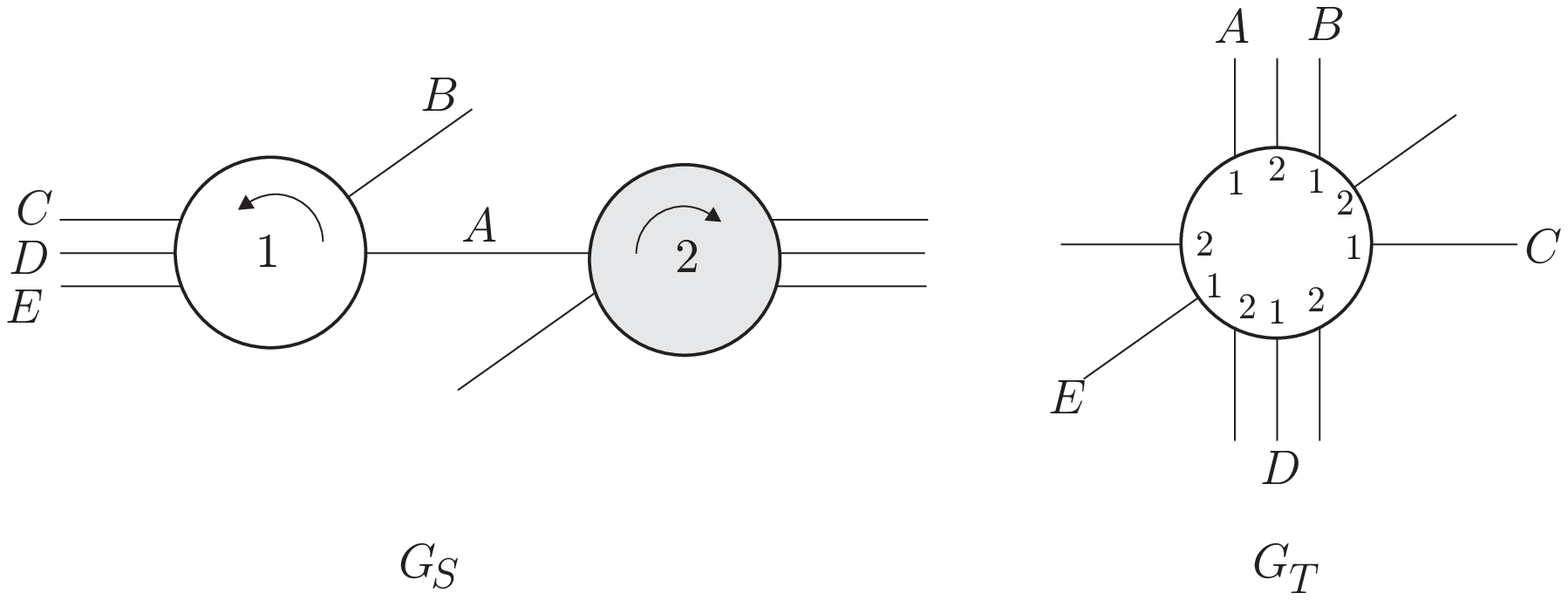}\vspace{-2mm}
\caption{}\label{fig:t1s2-jump1}
\end{center}
\end{figure}

Hence $G_T$ contains an $S$-cycle $\sigma_1$ consisting of edges $A, D$ whose face is $f_1$,
and a Scharlemann cycle $\sigma_2$ of length three with face $f_2$
consisting of edges $B, C, E$.
They lie on the same side of $\widehat{S}$.
Let us call this side the black side $\mathcal{B}$, and call the other the white side $\mathcal{W}$.
Let $H=V_\alpha\cap \mathcal{B}$.
Take $X=\widehat{S}\cup H\cup N(f_1\cup f_2)$ in $\mathcal{B}$.
Then $\partial X$ consists of the torus $\widehat{S}$ and the $2$-sphere.
For, $\partial f_1$ is non-separating on the genus two surface $F$ obtained from $\widehat{S}$ by tubing along $H$,
and $\partial f_2$ is non-separating on the torus obtained from $F$ by compressing along $f_1$.
Since $M(\alpha)$ is irreducible, its $2$-sphere bounds a ball in $\mathcal{B}$.
The situation in $\mathcal{W}$ is similar.
This means that $M(\alpha)$ is closed, and so $\partial M$ is a single torus.

The case where the jumping number is two is similar.
In fact, $G_S$ and $G_T$ are determined as shown in Figure \ref{fig:t1s2-jump2}, where
the correspondence of edges is indicated.
\end{proof}

\nocolon\begin{figure}[htb]
\begin{center}
\includegraphics*[scale=0.6]{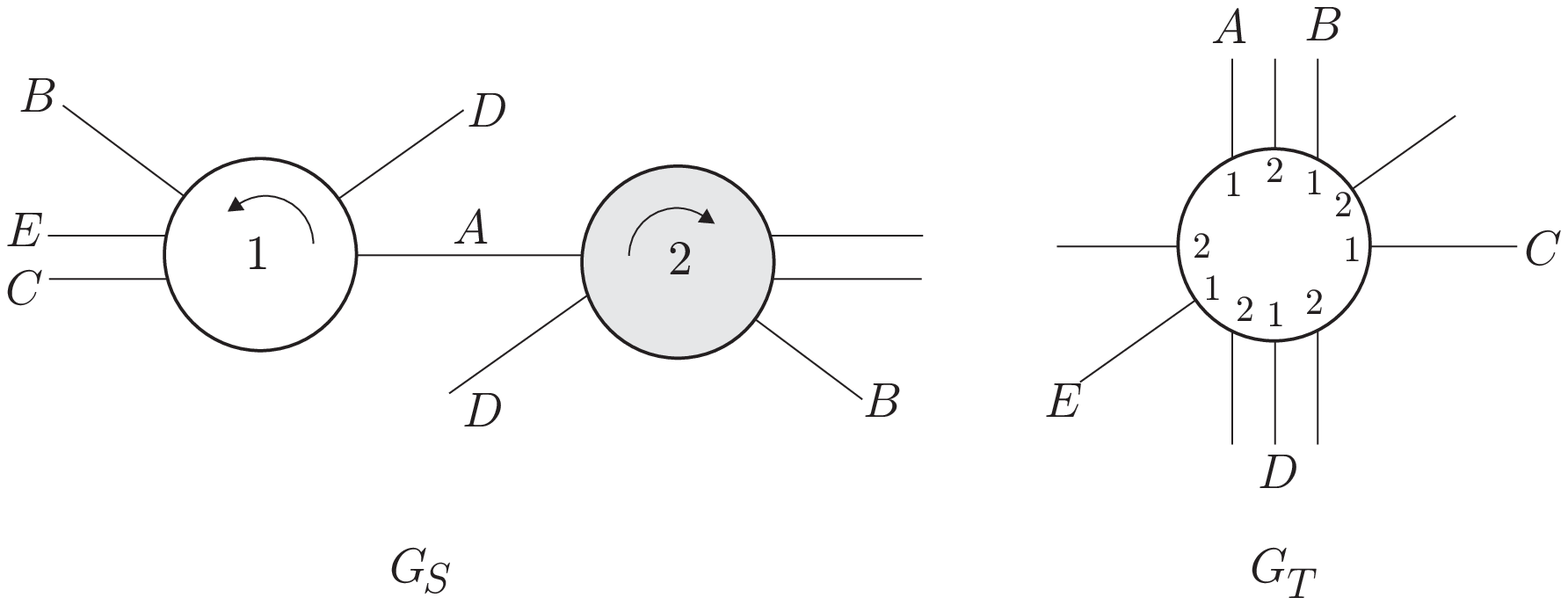}\vspace{-2mm}
\caption{}\label{fig:t1s2-jump2}
\end{center}
\end{figure}

Indeed, we can calculate $\pi_1M(\alpha)$ by Van Kampen's theorem.
Then if the jumping number is one, then $\pi_1M(\alpha)=\mathbb{Z}_5$, which contradicts that $M(\alpha)$ is toroidal.

\section{The case where $s\ge 4, t=2$}\label{sec:s>4t2}

\begin{figure}[htb]
\begin{center}
\includegraphics*[scale=0.3]{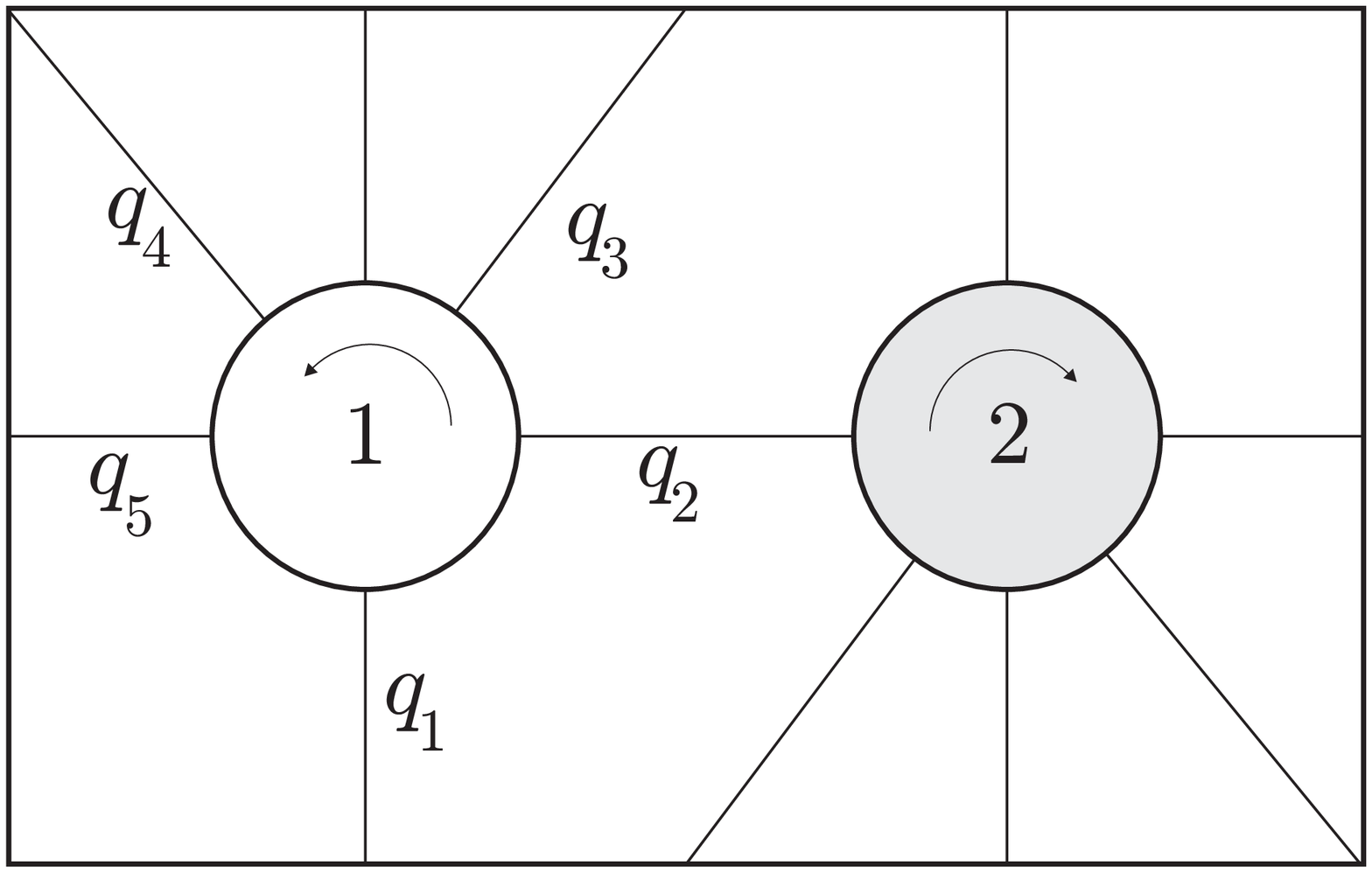}
\caption{}\label{fig:t2red}
\end{center}
\end{figure}

The reduced graph $\overline{G}_T$ is a subgraph of the graph as shown in Figure \ref{fig:t2red}.
Here, $q_i$ denotes the weight of edge.
As in \cite{Go2}, we say $G_T\cong G(q_1,q_2,q_3,q_4,q_5)$.

\nocolon

\begin{lemma}
Two vertices of $G_T$ have opposite signs.
\end{lemma}

\begin{proof}
Assume not.
Each vertex of $\overline{G}_T$ has degree at most $6$.
Note all edges of $G_T$ are positive.
Then $5s\le 6(s/2+1)=3s+6$, giving $s\le 3$, a contradiction.
\end{proof}

Then the arguments from Lemmas 5.1 to 5.3 of \cite{T} work with exchanging the roles of $G_S$ and $G_T$ there.
In particular, $q_1=s/2$ or $s/2+1$.
In the proofs of Lemmas 5.2 and 5.3 of \cite{T}, we use the fact that the jumping number is one.
But the case where the jumping number is two is similar.

\begin{lemma}\label{lem:half}
The case $q_1=s/2$ is impossible.
\end{lemma}

\begin{proof}
If $q_1=s/2$ then $G_T\cong G(s/2,s,s,s,s)$.
Then the four families of parallel negative edges correspond to the same permutation $\sigma$,
which is an involution \cite[Lemma 5.3]{T}.
Thus any component of $G_S^+$ has two vertices and $8$ edges.
Then there are $4$ mutually parallel edges, and so we have two bigons lying in the same side of $\widehat{T}$.
If these bigons do not have the same pair of edge class labels, then $M(\beta)$ contains a Klein bottle by \cite[Lemma 5.2]{GL3}.
Hence those have the same pair of edge class labels, but this is impossible by Lemma \ref{lem:common}.
\end{proof}


Thus $q_1=s/2+1$, and furthermore, Lemmas 5.8, 5.9 and 5.10 of \cite{T} hold.
(Instead of Lemma 2.7(2) of \cite{T}, we use the assumption that $M(\alpha)$ contains no Klein bottle.)
Hence we have $s=4$.
But this is shown to be impossible.

\begin{lemma}
$s=4$ is impossible.
\end{lemma}

\begin{proof}
We use the labeling of $G_T$ as in \cite[Figure 10]{T} (with changing $t$ to $s$).
In $G_S$, $u_1$ and $u_4$ are incident to three loops, and $u_2$ and $u_3$ are incident to two loops.
In $G_T$, there are two $S$-cycles with label pair $\{2,3\}$.
The edges of them give four edges between $u_2$ and $u_3$ in $G_S$.
Then two endpoints with label $1$ of loops at $u_2$ cannot be successive among five occurrences of label $1$.
Hence the jumping number is two.

In $G_S$, there are two edges between $u_1$ and $u_3$, which belong to $C$ in $G_T$.
Hence they are not parallel in $G_S$ by Lemma \ref{lem:common}.
Then there are two bigons at $u_1$ and $u_3$ which lie on the same side of $\widehat{T}$.
By \cite[Lemma 5.2]{GL3}, they must have the same pair of edge class labels.
Let $e$ be the remaining loop among three loops at $u_1$, not in the bigon.
By Lemma \ref{lem:common}, $e$ belongs to $D$ in $G_T$.
Also, let $c$ be the edge connecting $u_1$ and $u_3$ with the same label as $e$ at $u_1$.
Then the endpoints of $c$ and $e$ are consecutive at $v_1$ among five occurrences of label $1$,
which contradicts that the jumping number is two.
\end{proof}

\section{The case where $s=2, t>2$}

If $\widehat{T}$ is separating in $M(\beta)$, then the argument of Section \ref{sec:s>4t2} works with 
exchanging the role between $G_S$ and $G_T$.
Hence we suppose that $\widehat{T}$ is non-separating throughout this section.
We use $p_i$ to denote the weight in $\overline{G}_P$, instead of $q_i$ in Figure \ref{fig:t2red}.
Notice that $p_1\le t/2$, otherwise $G_S$ contains an $S$-cycle.

\begin{lemma}
$p_1=0$.
\end{lemma}

\begin{proof}
Assume $p_1\ne 0$.
Then $G_S$ contains a positive edge, and hence not all vertices of $G_T$ have the same sign.
By Lemma \ref{lem:parallel-max2}(2), $p_i\le t$ for $i=2,3,4,5$.
Since $2p_1+p_2+p_3+p_4+p_5=5t$, we have $p_1=t/2$ and $p_i=t$ for any $i\ne 1$.
Then Lemma \ref{lem:half} and the argument before it lead us to the conclusion.
\end{proof}

Thus the edges of $G_S$ are divided into at most $4$ edge classes, and then
some class contains more than $t$ edges.
This implies that all vertices of $G_T$ have the same sign by Lemma \ref{lem:parallel-max2}.
Also any edge of $G_T$ is a $\{1,2\}$-edge, and any disk face of $G_T$ is a Scharlemann cycle.

\begin{lemma}\label{lem:bandw}
$G_T$ has a black Scharlemann cycle and a white Scharlemann cycle.
\end{lemma}

\begin{proof}
Since $G_S$ has $5t$ edges, some edge class contains more than $t$ edges.
The associated permutation to the family has a single orbit by \cite[Lemma 4.2]{Go2}.
In particular, these $t+1$ edges cut $\widehat{T}$ into a disk.
Thus all faces of $G_T$ are disks, which gives a conclusion immediately.
\end{proof}

We say that two (disk) faces $f_1, f_2$ of $G_T$ of the same color are \textit{isomorphic\/} if the cyclic sequences
of edge class labels, read around their boundaries in the same direction, are equal.

\begin{lemma}\label{lem:atmost2}
$\partial M$ consists of at most two tori.
\end{lemma}

\begin{proof}
First, we prove:

\begin{claim}\label{cl:atmost3}
$\partial M$ consists of at most three tori.
\end{claim}

\begin{proof}[Proof of Claim \ref{cl:atmost3}]
Recall that $\widehat{S}$ separates $M(\alpha)$ into $\mathcal{B}$ and $\mathcal{W}$.
Let $H=V_\alpha\cap \mathcal{B}$ and let $f$ be a black Scharlemann cycle in $G_T$.
Then take a neighborhood $N=N(\widehat{S}\cup H\cup f)$ in $\mathcal{B}$.
Thus $\partial N=\widehat{S}\cup S'$, where $S'$ is a torus.
Since $S'\cap V_\alpha=\emptyset$, and $M$ is irreducible and atoroidal, either $S'$ bounds a solid torus in $\mathcal{B}$ or
$S'$ is parallel to a component of $\partial M$.
This means that $\partial\mathcal{B}$ consists of at most two torus boundary components, and
similarly for $\mathcal{W}$.
\end{proof}

Suppose that $\partial M$ consists of exactly three tori.
This happens only when both $\mathcal{B}$ and $\mathcal{W}$ have two tori as their boundaries.
Then all black disk faces of $G_T$ are isomorphic, and so are all white disk faces of $G_T$
by the argument of the proof of \cite[Lemma 5.6]{GL3}.
Notice that $G_T$ has $5t$ edges, but $\overline{G}_T$ has at most $3t$ edges, as seen by an easy Euler characteristic calculation.
Hence $G_T$ has a bigon.
Thus we may assume that all black faces are bigons.
By \cite[Lemma 5.2]{GL3}, all black bigons have the same pair of edge class labels, $\{\lambda,\mu\}$, say.
(For, $M(\alpha)$ contains no Klein bottle.)
Since all faces of $G_T$ are disks as in the proof of Lemma \ref{lem:bandw},
the set of edge class labels of any white face is also $\{\lambda,\mu\}$.
In $G_S$, this means that all edges are divided into two classes $\lambda$ and $\mu$.
Thus either of them contains more than $2t$ edges.
By \cite[Corollary 5.5]{Go2}, $t=3$.
Then $G_T$ has $15$ edges.
But this is impossible, because all black faces are bigons.
\end{proof}

\section{The case where $s=t=2$}\label{sec:st2}

The reduced graphs $\overline{G}_S$ and $\overline{G}_T$ are subgraphs of the graph as shown in Figure \ref{fig:t2red}.
Recall that we use $p_i$ (resp.\ $q_i$) to denote the weight of edge in $\overline{G}_S$ (resp.\ $\overline{G}_T$).

\subsection{Two vertices of $G_T$ have the same sign}

Since all edges in $G_T$ are positive, all edges in $G_S$ are negative.
Thus the edges of $G_S$ are divided into four edge classes.
Also, any edge of $G_T$ is a $\{1,2\}$-edge, and any disk face of $G_T$ is a Scharlemann cycle.

\begin{lemma}
$G_T$ has a black Scharlemann cycle and a white Scharlemann cycle.
\end{lemma}

\begin{proof}
If some $q_i>2$, then $G_T$ contains a black bigon and a white bigon.
Hence we assume that $q_i\le 2$ for any $i$.
Since $2q_1+q_2+q_3+q_4+q_5=10$, $q_1\ne 0$.
If $q_1=1$, then $q_2=q_3=q_4=q_5=2$, giving the conclusion.
If $q_1=2$, then we can assume that $q_2+q_3=4$ and $q_4+q_5=2$ by symmetry.
Then $q_2=q_3=2$, giving the conclusion.
\end{proof}

By the same argument in the proof of Claim \ref{cl:atmost3}, $\partial M$ consists of at most three tori.

\begin{lemma}
$\partial M$ consists of at most two tori.
\end{lemma}

\begin{proof}
If not, then as in the proof of Lemma \ref{lem:atmost2},
all black disk faces of $G_T$ are isomorphic, and so are all white disk faces of $G_T$.
If $q_i\ge 3$ for some $i$, then all disk faces of $G_T$ would be bigons, which is impossible.
Hence $q_i\le 2$ for any $i$.
In particular, $q_1>0$.
If $q_1=1$, then $q_2=q_3=q_4=q_5=2$, which
contradicts that all disk faces of the same color are isomorphic.
If $q_1=2$, then we may assume that $q_2=q_3=2$ and $q_4+q_5=2$ by symmetry.
But any case where $(q_4,q_5)=(2,0), (1,1)$ gives a contradiction similarly.
\end{proof}

\subsection{Two vertices of $G_T$ have distinct signs}

We will show that there is only one possible pair for $\{G_S,G_T\}$.
Lemmas 6.1 and 6.2 of \cite{T} hold here (the jumping number two case is similar in the proof of Lemma 6.2 of \cite{T}),
and hence $p_1=2$ or $3$.

\begin{figure}[htb]
\begin{center}
\includegraphics*[scale=0.6]{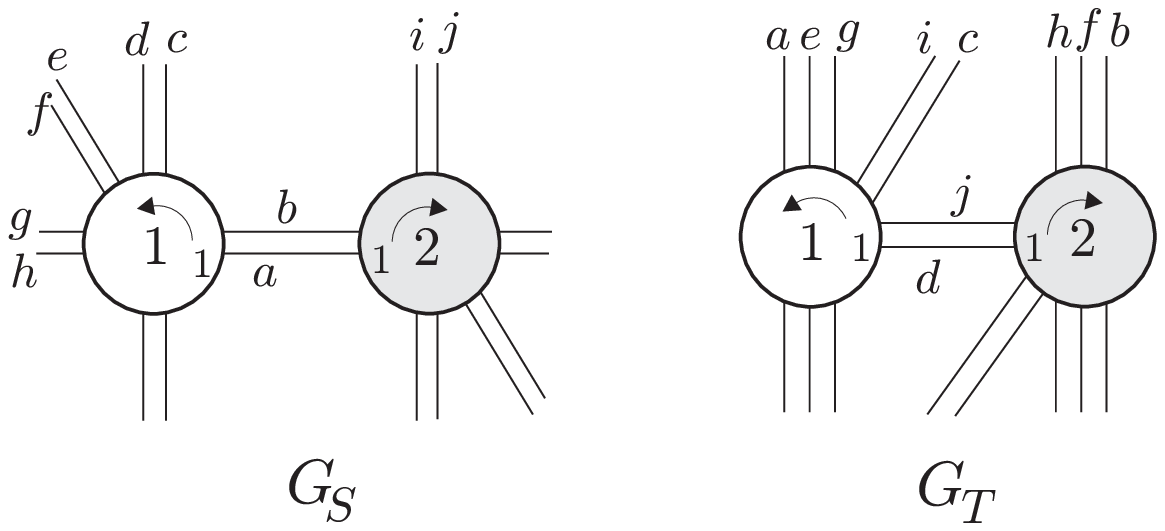}\vspace{-2mm}
\caption{}\label{fig:s2t2final}
\end{center}
\end{figure}

\begin{lemma}
If $p_1=2$, then the graphs are as shown in Figure \textup{\ref{fig:s2t2final}}, where the jumping number is two.
\end{lemma}

\nocolon

\begin{proof}
As in the proof of \cite[Lemma 6.3]{T},
there is only one possibility for $G_T$ as shown in \cite[Figure 16(4)]{T}.
In fact, if the jumping number is one, then this is also eliminated as shown there.
\end{proof}

\begin{lemma}
If $p_1=3$, then the graphs are the same as in Figure \textup{\ref{fig:s2t2final}} with exchanging $G_S$ and $G_T$.
\end{lemma}

\begin{proof}
We may assume that $(p_2+p_3,p_4+p_5)=(4,0)$ or $(2,2)$ by symmetry.
In the latter case, there are three possibilities for $G_S$ as in the proof of \cite[Lemma 6.4]{T}, and
all are impossible.
Thus $(p_2+p_3,p_4+p_5)=(4,0)$, giving $p_2=p_3=2$.
Hence $q_1=2$, and so we can assume that $(q_2+q_3,q_4+q_5)=(6,0)$ or $(4,2)$ by symmetry.
Then $(6,0)$ contradicts Lemma \ref{lem:common}.
By using the parity rule,
it is easy to see that $G_T$ is as in Figure \ref{fig:s2t2final} (with exchanging $G_S$ and $G_T$).
\end{proof}

\begin{lemma}
If the graphs are as in Figure \textup{\ref{fig:s2t2final}}, then $\partial M$ consists of at most two tori.
\end{lemma}

\begin{proof}
We may use the notation of Figure \ref{fig:s2t2final}.
Then we can assume that $G_T$ contains $4$ black bigons and two white bigons and two white $3$-gons.
As in the proof of Claim \ref{cl:atmost3}, the black side $\mathcal{B}$ of $\widehat{S}$ in $M(\alpha)$
has $\widehat{S}$ and at most one torus as its boundary.
On the other hand, the white side $\mathcal{W}$ has a single torus $\widehat{S}$ as its boundary,
because a torus obtained from $\widehat{S}\cup (V_\alpha\cap \mathcal{W})$ by attaching a bigon, will be compressed
by a $3$-gon.
Hence $\partial M$ consists of at most two tori.
\end{proof}

\section{Klein bottle}\label{sec:klein}

In the rest of the paper, we will treat the case where either $M(\alpha)$ or $M(\beta)$ contains
a Klein bottle.

Suppose that $M(\alpha)$ contains a Klein bottle $\widehat{P}$ such that $\widehat{P}\cap V_\alpha$
consists of $p$ meridian disks $u_1,u_2,\dots,u_p$, numbered successively, of $V_\alpha$,
and that $p$ is minimal among all Klein bottles in $M(\alpha)$.
Let $P=\widehat{P}\cap M$.  Since $M$ is hyperbolic, $p>0$.
Remark that we do not assume that $M(\alpha)$ is toroidal.

Now, we suppose $\Delta=\Delta(\alpha,\beta)\ge 6$.
Then notice that both of $M(\alpha)$ and $M(\beta)$ are irreducible by \cite{LOT,O,O2,W}.
Let $S'=\partial N(\widehat{P})$.
If $S'$ is boundary parallel in $M(\alpha)$, then $M(\alpha)=N(\widehat{P})$, and hence
$\partial M$ consists of two tori.
If $M(\alpha)$ is also toroidal, then $\partial M$ is at most two tori by \cite{Go2}.
Hence we assume $S'$ is compressible in $M(\alpha)$.
But this implies that $S'$ bounds a solid torus by the irreducibility of $M(\alpha)$, and so $\partial M$ is a single torus.
Therefore we assume that $\Delta=5$ in the rest of the paper.
Then both of $M(\alpha)$ and $M(\beta)$ are irreducible.

\begin{lemma}\label{lem:incom}
$P$ is incompressible and boundary-incompressible in $M$.
\end{lemma}

\begin{proof}
Suppose that $P$ is compressible in $M$.
Let $D$ be a disk in $M$
such that $D\cap P=\partial D$ and $\partial D$ does not bound a disk on $P$.
Note that $\partial D$ is orientation-preserving on $P$.

If $\partial D$ is non-separating on $\widehat{P}$, then
we get a non-separating $2$-sphere in $M(\alpha)$ by compressing $\widehat{P}$ along $D$.
This contradicts the irreducibility of $M(\alpha)$.
If $\partial D$ bounds a disk on $\widehat{P}$, then we replace the disk with $D$, and
get a new Klein bottle in $M(\alpha)$ with fewer intersections with $V_\alpha$ than
$\widehat{P}$.
This contradicts the choice of $\widehat{P}$.
Thus $\partial D$ is essential and separating on $\widehat{P}$.
Compressing $\widehat{P}$ along $D$ gives two disjoint projective planes in $M(\alpha)$.
Since $M(\alpha)$ is irreducible, this is also impossible.
Thus we have shown that $P$ is incompressible.

Next, let $E$ be a disk in $M$
such that $E\cap P=\partial E\cap P$, $\partial E=a\cup b$, where
$a\subset P$ is an essential (i.e., not boundary-parallel) arc in $P$
and $b\subset \partial M$.
If $a$ joins distinct components of $\partial P$, then a compressing disk for $P$ is
obtained from two parallel copies of $E$ and the disk obtained by removing a neighborhood
of $b$ from the annulus in $\partial M$ cobounded by those components of $\partial P$
meeting $a$. 
Hence $\partial a$ is contained in the same component $\partial_1P$, say, of $\partial P$.
If $p>1$, then $b$ bounds a disk $D'$ in $\partial M$ together with a subarc of $\partial_1P$.
Then $E\cup D'$ gives a compressing disk for $P$ in $M$.
Therefore $p=1$. Then we can move the core of $V_\alpha$ onto an
orientation-reversing loop in $\widehat P$ by using $E$.
This implies that $M$ contains a properly embedded M\"{o}bius band,
which contradicts the fact that $M$ is hyperbolic. 
\end{proof}

Thus we can define two graphs $G_P$ on $\widehat{P}$ and $G_T$ on $\widehat{T}$
from the arcs in $P\cap T$ as in Section \ref{sec:pre}.
We can label each endpoint of edges of these graphs as before.
Note that neither $G_P$ nor $G_T$ has a trivial loop.
Lemma \ref{lem:jumping} holds without any change.

Since $\widehat{P}$ is non-orientable, we cannot give a sign to a vertex of $G_P$ as in $G_T$.
Hence assign an orientation to each vertex of $G_P$ as a meridian disk of $V_\alpha$.
That is, all vertices of $G_P$ determine the same homology class in $H_2(V_\alpha,\partial V_\alpha)$.
By using this, we give a sign to each edge of $G_P$ as follows.

Let $e$ be an edge of $G_P$.
Assume that $e$ is a loop based at $u$.
Then $e$ is \textit{positive\/} if a regular neighborhood $N(u\,\cup\,e)$ on $\widehat{P}$ is an annulus, \textit{negative\/} otherwise.
Assume that $e$ connects distinct vertices $u_i$ and $u_j$.
Then $N(u_i\,\cup\, e\,\cup\,u_j)$ is a disk.
Then $e$ is \textit{positive\/} if we can give an orientation to the disk $N(u_i\,\cup\, e\,\cup\, u_j)$ so that
the induced orientations on $u_i$ and $u_j$ are compatible with those of $u_i$ and $u_j$ simultaneously.
Otherwise, $e$ is \textit{negative\/}.
Then the parity rule (Lemma \ref{lem:common}(1)) still holds without change.
In fact, the above definition works for $G_T$, and so this is a natural generalization of
the usual parity rule.
Also, Lemma \ref{lem:common}(2) is true.


\begin{lemma}\label{lem:parallel-max-kb}
$G_P$ satisfies the following\textup{:}
\begin{itemize}
\item[\rm(1)] If $t\ge 3$, then any family of parallel positive edges contains at most $t/2+2$ edges.
Moreover, if it contains $t/2+2$ edges, then $t\equiv 0 \pmod{4}$, and, up to relabelling of vertices of $G_T$,
it contains $\{1,2\}$ $S$-cycle and $\{t/2,t/2+1\}$ $S$-cycle.
\item[\rm(2)] Either all the vertices of $G_T$ have the same vertex, or
any family of negative edges contains at most $t$ edges.
In particular, if $G_P$ contains a positive edge, any family of negative edges contains at most $t$ edges.
\end{itemize}
\end{lemma}

\begin{proof}
(1) is \cite[Lemma 1.4 and Corollary 1.8]{W}.
(2) is the same as Lemma \ref{lem:parallel-max2}(2).
\end{proof}

If $p\ge 3$, a \textit{generalized $S$-cycle\/} in $G_T$ is the triplet of mutually parallel positive
edges $e_{-1},e_0,e_1$, where $e_{-1}$ and $e_1$ have the same label pair $\{i-1,i+1\}$, and $e_0$
is a level $i$-edge for some $i$.

\begin{lemma}\label{lem:kb-S}
$G_T$ has neither a Scharlemann cycle nor a generalized $S$-cycle.
\end{lemma}

\begin{proof}
For a Scharlemann cycle, see \cite[Lemma 3.2]{T0}.
(It treats the case of $S$-cycles, but the argument works for general case.)
\end{proof}

\begin{lemma}\label{lem:parallel-max3}
Assume $p\ge 2$.  Then $G_T$ satisfies the following\textup{:}
\begin{itemize}
\item[\rm(1)] Any family of parallel positive edges contains at most $p/2+1$ edges.
Moreover, if it contains $p/2+1$ edges, then the first and last edge are level.
\item[\rm(2)] Any family of parallel negative edges contains at most $p$ edges.
\end{itemize}
\end{lemma}

\begin{proof}
(1) Assume that $G_T$ contains a family $A$ of mutually parallel positive edges
which connect $v_i$ and $v_j$ (possibly, $i=j$), and that
$A$ contains more than $p/2+1$ edges.

When $p=2$, no edge of $A$ is level.
Otherwise, there would be a pair of edges which are parallel in both graphs.
But this is impossible by Lemma \ref{lem:common}(2).
Hence $A$ contains an $S$-cycle, a contradiction by Lemma \ref{lem:kb-S}.

Suppose $p>2$.
Note that some label appears at both $v_i$ and $v_j$.
If $A$ contains no level edge, then $A$ contains
an $S$-cycle.
This is impossible by Lemma \ref{lem:kb-S}.
Hence $A$ must contain a level edge.
Moreover, a level edge is the first or last edge of $A$.
Otherwise, $A$ contains a generalized $S$-cycle, which is also impossible by Lemma \ref{lem:kb-S}.
We may assume that the first edge of $A$ is level.
Then $A$ contains an $S$-cycle if $p$ is odd.
If $p$ is even, the second to last edge is level, and hence there is a generalized
$S$-cycle.
The second assertion is easy to see.

(2)\qua Assume that $G_T$ contains $p+1$ parallel negative edges, connecting $v_i$ and $v_j$.
Consider the associated permutation $\sigma$ to these edges as follows.
Let $a_1,a_2,\dots,a_p,b_1$ be the edges labelled successively.
We may assume that $a_k$ has label $k$ at $v_i$, label $\sigma(k)$ at $v_j$.
Let $\theta$ be the orbit of $\sigma$ containing $1$,
and let $C_\theta$ be the cycle in $G_P$ corresponding to $\theta$.
Then $C_\theta$ does not bound a disk in $\widehat{P}$ by \cite[Lemma 2.3]{Go2}.
Note that there are two possibilities for $C_\theta$, that is,
separating or non-separating in $\widehat{P}$, since $C_\theta$ is orientation-preserving
in $\widehat{P}$.
Consider the edge $b_1$.
Since $b_1$ is positive in $G_P$,
either $b_1$ is parallel to $a_1$ in $G_P$, or,
the cycle consisting of the edges $a_2,\dots,a_p,b_1$ bounds a disk in $\widehat{P}$.
But the former contradicts Lemma \ref{lem:common}(2), and the latter is impossible by \cite[Lemma 2.3]{Go2} again.
\end{proof}

In this section, we treat the case that $G_P$ or $G_T$ has a single vertex.

\begin{proposition}\label{prop:t1}
If $t=1$, then $\partial M$ is a single torus.
\end{proposition}

\begin{proof}
Suppose $t=1$.
If $p=1$, then $G_P$ has a single vertex with degree $5$, which is impossible.
Recall that the edges of $G_T$ are divided into at most three edge classes as in Section \ref{sec:t1}.
Also, at most $p/2+1$ edges can be parallel in $G_T$.
If $p\ge 3$, then $6(p/2+1)\ge 5p$ gives $p=3$.
But then at most two edges can be parallel in $G_T$, giving $6\cdot 2\ge 15$, a contradiction.
Thus $p=2$, and hence $G_T\cong G(2,2,1)$. (Recall the notation in Section \ref{sec:t1}.)
Then $G_S$ is uniquely determined, and the correspondence between the edges of $G_P$ and $G_T$ is shown in Figure \ref{fig:kbt1}.
Here, the jumping number must be one, and
two end circles of the cylinder are identified through a suitable involution
to form the Klein bottle $\widehat{P}$.
Note that the edge connecting two vertices with labels $1$ and $2$ is negative in $G_P$.

\nocolon\begin{figure}[htb]
\begin{center}
\includegraphics*[scale=0.33]{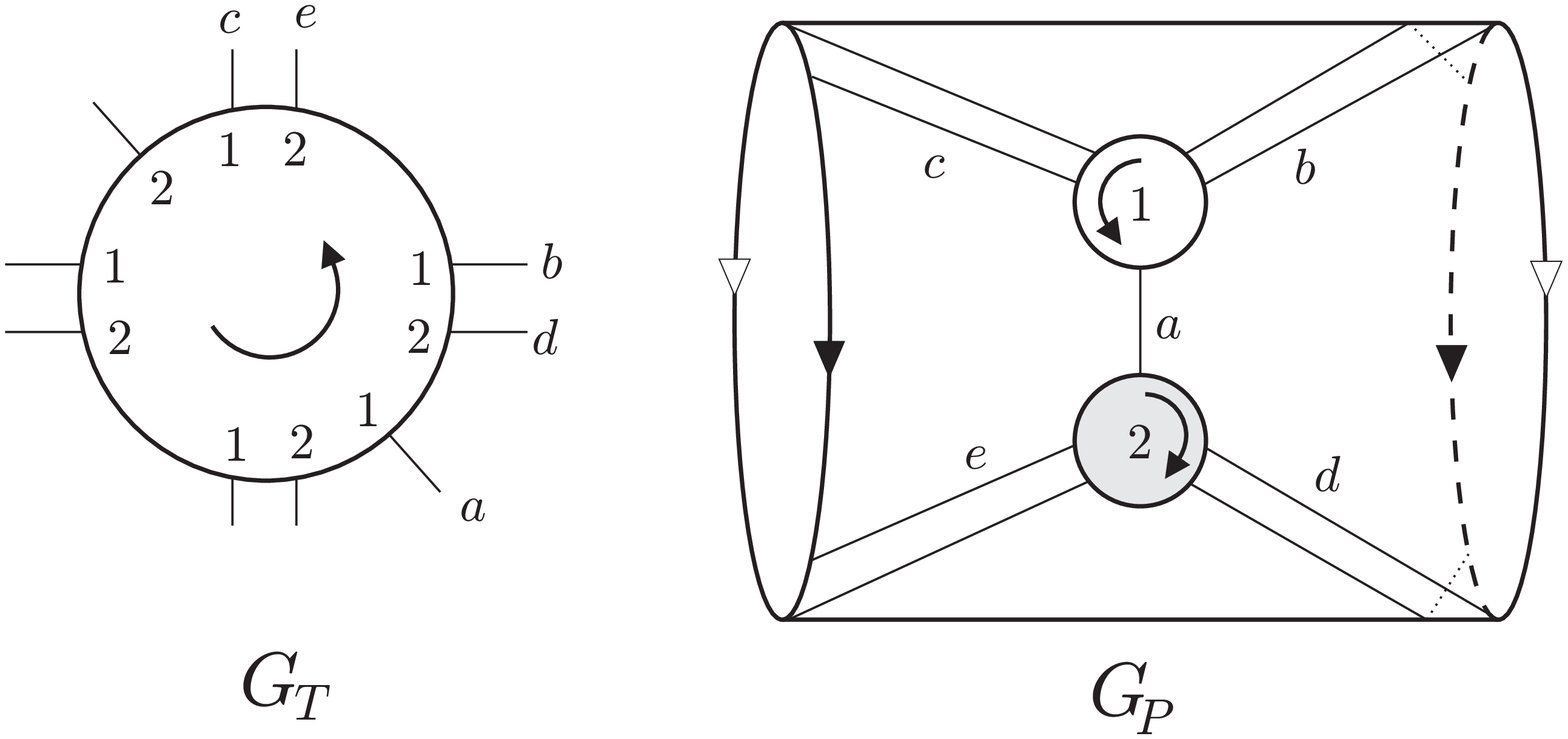}\vspace{-2mm}
\caption{}\label{fig:kbt1}
\end{center}
\end{figure}

Let $N(\widehat{P})$ be a regular neighborhood of $\widehat{P}$ in $M(\alpha)$.
Then $N(\widehat{P})$ is the twisted $I$-bundle over $\widehat{P}$, and
$\partial N(\widehat{P})$ is a torus.
Let us write $M(\alpha)=N(\widehat{P})\cup W$.
Then $T\cap W$ consists of two bigons and two $3$-gons.
Also, $V_\alpha\cap W$ consists of two $1$-handles $H_1, H_2$.
Let $F$ be the genus three closed surface obtained from $\partial W$ by performing
surgery along $H_1$ and $H_2$.
Then attaching a bigon and two $3$-gons to $F$ yield the $2$-sphere.
Since $M(\alpha)$ is irreducible, $M(\alpha)$ must be closed.
The result immediately follows from this.
\end{proof}

\begin{figure}[htb]
\begin{center}
\includegraphics*[scale=0.41]{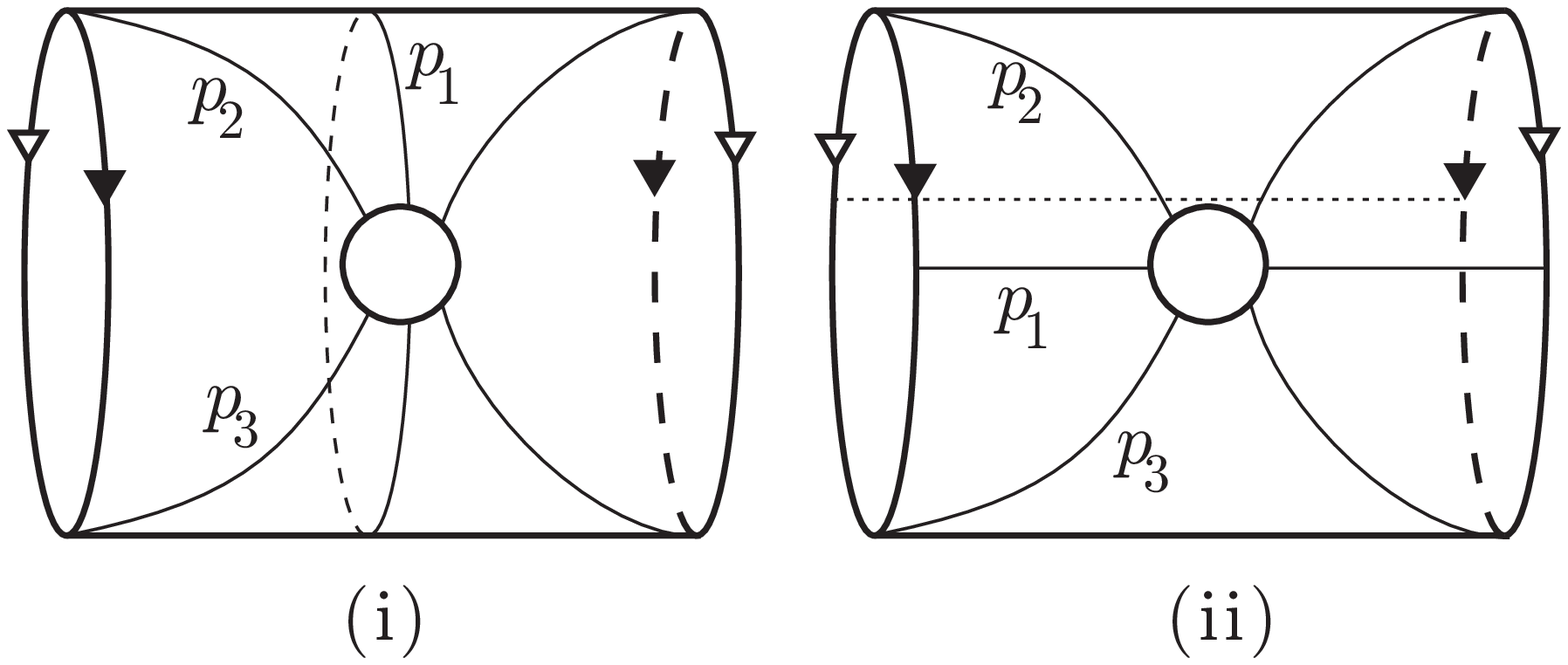}\vspace{-2mm}
\caption{}\label{fig:p1}
\end{center}
\end{figure}

\begin{lemma}
If $p=1$, then $\overline{G}_P$ is a subgraph of either graph shown in Figure \textup{\ref{fig:p1}}.
\end{lemma}

\nocolon

\begin{proof}
An orientation-preserving loop on a Klein bottle is non-separating or separating.
Also, there are only two classes of orientation-reversing loops.
The result follows immediately. 
(See \cite[Lemma 2.1]{P}.)
\end{proof}

Thus we say $G_P\cong H(p_1,p_2,p_3)$ for (i), or $H'(p_1,p_2,p_3)$ for (ii), where
$p_1$ denotes the weight of the positive loop, and the others denote the weight of negative loops in each class.
Clearly, $H(p_1,p_2,p_3)\cong H(p_1,p_3,p_2)$, $H'(p_1,p_2,p_3)\cong H'(p_1,p_3,p_2)$ and $H(0,p_2,p_3)\cong H'(0,p_2,p_3)$.
Also, $2(p_1+p_2+p_3)=5t$ implies that $t$ is even.

\begin{proposition}\label{prop:p1t2}
If $p=1$ and $t=2$, then $\partial M$ consists of at most two tori.
\end{proposition}

\begin{proof}
First, we claim $p_1\ne 0$.
If $p_1=0$, then $G_P\cong H(0,5,0)$, $H(0,4,1)$ or $H(0,3,2)$.
For $H(0,5,0)$, $G_T$ contains $5$ edges connecting $v_1$ and $v_2$.
Since there are at most $4$ edge classes, this contradicts Lemma \ref{lem:common}(2).
For $H(0,4,1)$, $G_T$ has two loops at each vertex, which must be parallel.
So, this contradicts Lemma \ref{lem:common}(2) again.
For $H(0,3,2)$, $G_T\cong G(1,1,1,1,0)$ by using Lemma \ref{lem:common}(2).
Then a jumping number argument eliminates this as follows.
By examining the endpoints of a loop at $v_1$, we see that the jumping number is two.
Let $a$ and $b$ be the edges connecting $v_1$ and $v_2$ such that
their end points at $v_1$ are consecutive.
Then they are parallel in $G_P$ and adjacent.
(In fact, they belong to the family of $3$ mutually parallel negative edges of $G_P$.)
By Lemma \ref{lem:jumping}, the endpoints of $a$ and $b$ with label $1$ are not
consecutive at $u_1$ among five occurrences of label $1$.
Then their endpoints with label $2$ are consecutive among five occurrences of label $2$ at $u_1$.
But $a$ and $b$ are consecutive at $v_2$ also, which contradicts Lemma \ref{lem:jumping}.

Notice that $1\le p_1\le 5$.
If $p_1=5$, then we have a pair of edges which are parallel in both graphs, a contradiction.
In the following, we consider all possibilities for $G_P$.

Seven cases $H(4,1,0)$, $H(3,2,0)$, $H(2,3,0)$, $H(2,2,1)$, $H(1,4,0)$, $H(1,2,2)$, $H'(1,3,1)$ are impossible by the parity rule.
For the four cases $H'(4,1,0)$, $H'(3,1,1)$, $H'(2,3,0)$, $H'(2,2,1)$, $G_P$ contains an $S$-cycle.
Hence $\widehat{T}$ is separating, and so the faces of $G_P$ can be colored by two colors in such a way that
two sides of an edge have distinct colors.
This fact eliminates these four cases.
For $H(1,3,1)$ and $H'(1,4,0)$, $G_T$ contains two loops which are parallel in both graphs.

For $H'(1,2,2)$, $G_T\cong G(2,1,0,0,0)$.
At $v_1$, there is no correct arrangement of edges to satisfy Lemma \ref{lem:jumping}.
For $H(3,1,1)$, $G_T\cong G(1,1,1,1,0)$.
As in the proof of Proposition \ref{prop:t1}, 
let $M(\alpha)=N(\widehat{P})\cup W$.
Then $T\cap W$ contains two $3$-gons.
Attaching these $3$-gons to $N(\partial W\cup (V_\alpha\cap W))$ yields a $2$-sphere.
Since $M(\alpha)$ is irreducible, $M(\alpha)$ is closed.
Thus $\partial M$ is a single torus.  
Finally, for $H'(3,2,0)$, $G_T\cong G(1,1,1,1,0)$ again.
Take one $3$-gon in $T\cap W$.
Attaching it to $N(\partial W\cup (V_\alpha\cap W))$ yields a torus $S'$, missing $V_\alpha$.
Thus $S'$ is boundary parallel or compressible.
In the former, $\partial M$ consists of two tori.
In the latter, either $S'$ bounds a solid torus in $W$,
which implies that $\partial M$ is a single torus,
or $S'$ is contained in a $3$-ball in $M(\alpha)$, which implies
that $S'$ bounds a knot exterior $X$.
Since a Klein bottle cannot lie in a knot exterior, $X$ lies in $W$.
In any case, $\partial M$ is a single torus.
\end{proof}

\begin{proposition}
If $p=1$ and $t>2$, then $\partial M$ consists of at most two tori.
\end{proposition}

\begin{proof}
By Lemma \ref{lem:parallel-max-kb}, $p_1\le t/2+2$.
Hence $p_2+p_3=5t/2-p_1\ge 2t-2$.
Then an Euler characteristic calculation shows that $G_T^+$ has a disk face $D$.
Let us write $M(\alpha)=N(\widehat{P})\cup W$.
Then $\partial N(\partial W\cup (V_\alpha\cap W)\cup D)$ consists of two tori, since $\partial D$ runs on the $1$-handle $V_\alpha\cap W$ in the same
direction.
This implies that $\partial M$ consists of at most two tori as in the last paragraph of the proof of Proposition \ref{prop:p1t2}.
\end{proof}

\section{Klein bottle; the case $t=2$}

By Section \ref{sec:klein}, we may assume $p\ge 2$.

\begin{lemma}
Two vertices of $G_T$ have opposite signs.
\end{lemma}

\begin{proof}
Assume not.  Then $q_i\le p/2+1$ for any $i$.
Thus $5p\le 6(p/2+1)=3p+6$ gives $p\le 3$.
If $p=3$, then $q_i\le 2$, and so $15=5p\le 12$, a contradiction.
Assume $p=2$.
Since $q_i\le 2$ for any $i$, $q_1=1$ or $2$.
Hence $G_T\cong G(1,2,2,2,2)$ or $G(2,q_2,q_3,q_4,q_5)$ with $q_2+q_3=q_4+q_5=3$.

For $G(1,2,2,2,2)$, the labels of $G_T$ are determined, up to exchange of $1$ and $2$, and then $G_P$ is uniquely determined.
See Figure \ref{fig:kb-t2-p2}.
Consider the edges $a$, $b$ and $c$ as there.
The endpoints of $a$ and $c$ are consecutive at $v_1$, but those of $b$ and $c$ are not consecutive at $v_2$,
among the five occurrences of label $1$.
Any location of $c$ contradicts Lemma \ref{lem:jumping} at $u_1$.

\nocolon\begin{figure}[htb]
\begin{center}
\includegraphics*[scale=0.29]{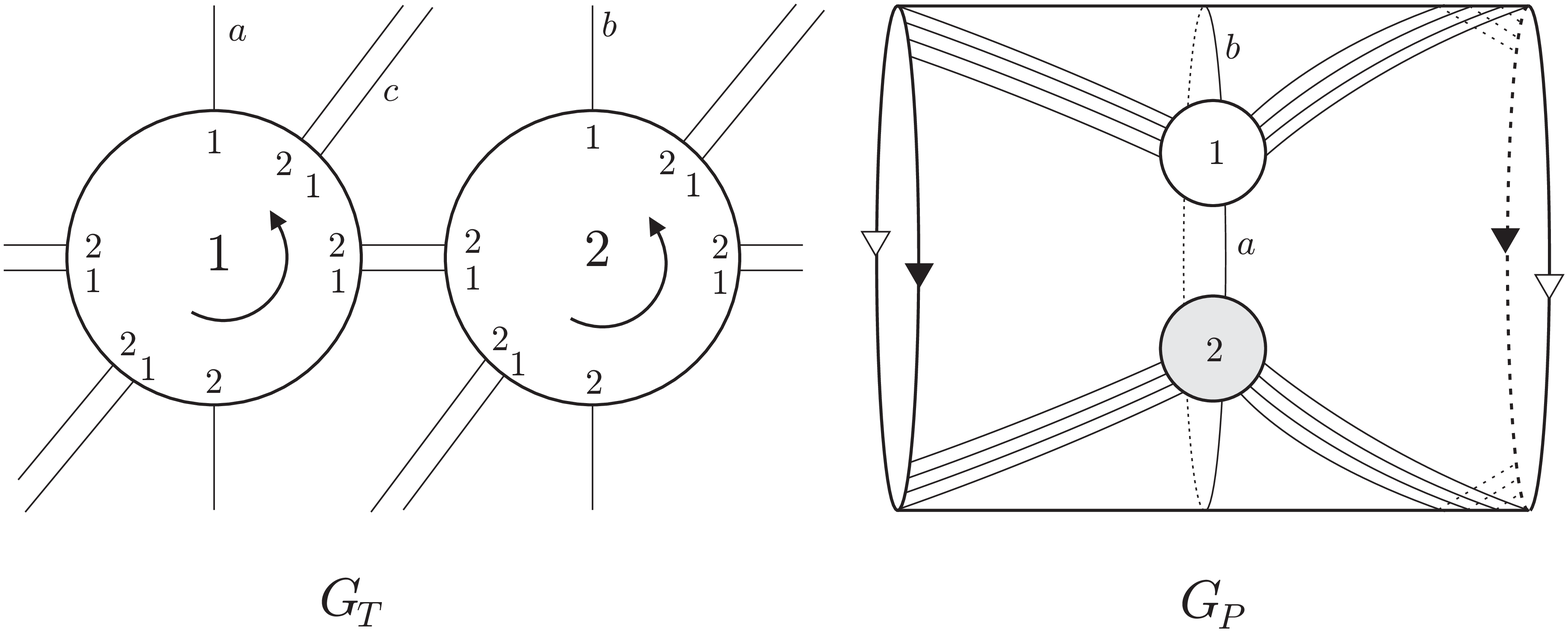}\vspace{-2mm}
\caption{}\label{fig:kb-t2-p2}
\end{center}
\end{figure}

Suppose $G_T\cong G(2,q_2,q_3,q_4,q_5)$ with $q_2+q_3=q_4+q_5=3$.
Then $G_T\cong G(2,2,1,2,1)$ or $G(2,2,1,1,2)$.
In any case, each vertex of $G_P$ is incident to $4$ negative loops, where are parallel.
But two of them are level, and the others are not level, a contradiction.
\end{proof}

\begin{lemma}
If $p$ is even, then $q_1=p/2$ or $p/2+1$.
If $p$ is odd, then $q_1=(p+1)/2$.
\end{lemma}

\begin{proof}
By Lemma \ref{lem:parallel-max3}, $q_1\le p/2+1$ and $q_i\le p$ for $i\ne 1$.
Since $2q_1+q_2+q_3+q_4+q_5=5p$, we have $q_1\ge p/2$, giving the conclusion.
\end{proof}

We consider three cases.

\subsection{$q_1=p/2$}

Then $G_T\cong G(p/2,p,p,p,p)$.
Let $A_i$ be the family of parallel negative edges of weight $q_i$ for $i=2,3,4,5$.
Then they associate to the same permutation $\sigma$.

\begin{lemma}
$\sigma$ is not the identity.
\end{lemma}

\begin{figure}[htb]
\begin{center}
\includegraphics*[scale=0.36]{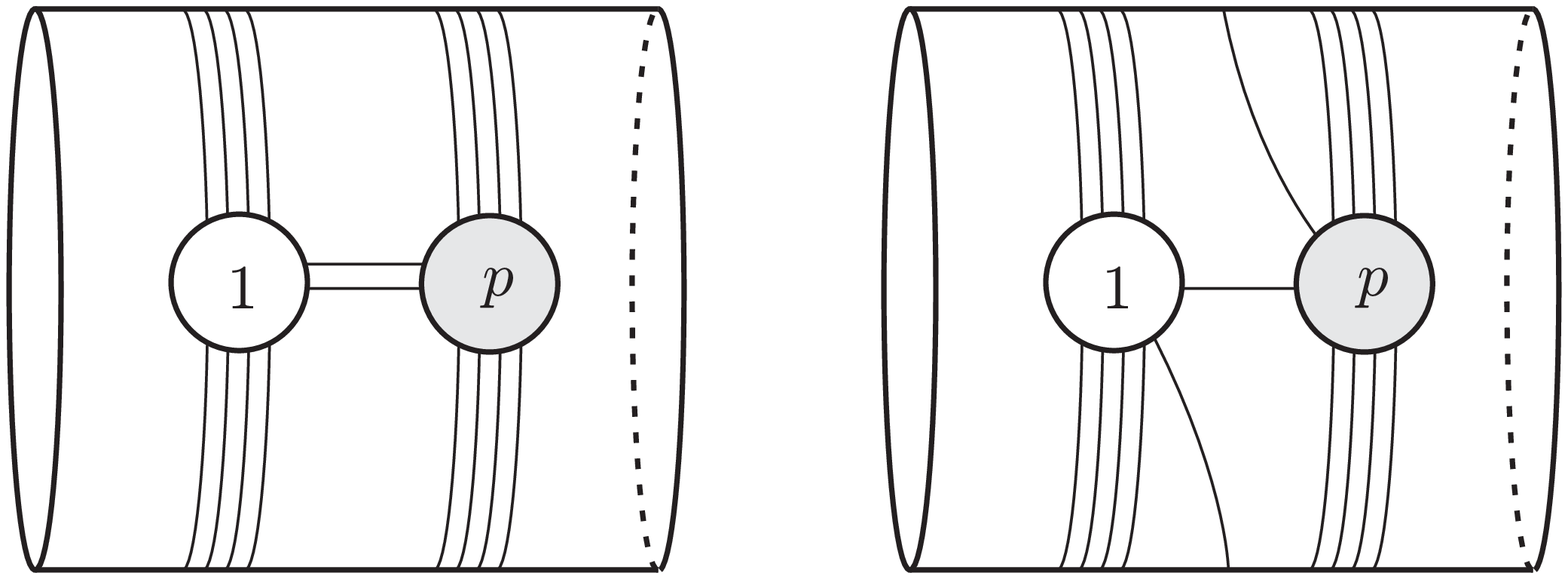}
\caption{}\label{fig:kbt2}
\end{center}
\end{figure}

\begin{proof}
If $\sigma$ is the identity, then each family $A_i$ contains a $\{1,1\}$-edge and $\{p,p\}$-edge.
Let $G(1,p)$ be the subgraph of $G_P$ spanned by $u_1$ and $u_p$.
Then $G(1,p)$ has the form as in Figure \ref{fig:kbt2}.
But a jumping number argument gives a contradiction.
\end{proof}

\nocolon

\begin{lemma}
If $p=2$, then $\partial M$ consists of at most two tori.
\end{lemma}

\begin{proof}
Let us write $M(\alpha)=N(\widehat{P})\cup W$ again.
Then $T\cap W$ consists of four bigons and four $3$-gons.
Let us choose a bigon $D_1$ and a $3$-gon $D_2$.
It is easy to see that if $X=N(\partial W\cup (V_\alpha\cap W)\cup D_1\cup D_2)$
then $\partial X=\partial W\cup S'$, where $S'$ is a torus missing $V_\alpha$.
The result follows from this as in the proof of Proposition \ref{prop:p1t2}.
\end{proof}

Hence we assume $p\ge 3$ hereafter.

\begin{lemma}\label{lem:involution}
If $\sigma$ is not the identity, then $\sigma^2$ is the identity.  In particular, each orbit of $\sigma$ has length two.
\end{lemma}

\begin{proof}
The proof of \cite[Lemma 5.3]{T} works here.
\end{proof}

\begin{lemma}
$q_1=p/2$ is impossible.
\end{lemma}

\begin{proof}
We may assume that the edges of $A_1$ have labels $1,2,\dots,p$ at $u_1$.
We follow the argument of \cite[Lemma 5.4]{T}. 
Then the component $H$ of $G_P$ containing $G(1,p/2+1)$ and $G(p/2,p/2+1)$ has
the form as in Figure 11 of \cite{T}.
(Here, we do not need the assumption $p>4$.)
But a jumping number argument eliminates this configuration
(even for the case that the jumping number is two).
\end{proof}

\subsection{$q_1=p/2+1$}

Since $G_T$ cannot contain a Scharlemann cycle,
$G_T\cong G(p/2+1,p,p-1,p,p-1)$ or $G(p/2+1,p,p-1,p-1,p)$.
Notice that the first and last edges of the positive loops at each vertex of $G_T$ are level.
We may assume that the edges of $A_1$ have labels $1,2,\dots,p$ at $u_1$.
Let $\sigma$ be the associated permutation to $A_1$.
Then $\sigma(i)\equiv i-1$ or $i+p/2-1 \pmod{p}$, since $p/2$ and $p$ are the only labels
of positive level edges in $G_T$.

\begin{lemma}
If $p>2$, then $\sigma(i)\equiv i+p/2-1 \pmod{p}$.
\end{lemma}

\begin{proof}
Assume $\sigma(i)\equiv i-1\pmod{p}$.
Then the edges of $A_1$ form an essential cycle $C$ through all vertices, which is separating or non-separating on $\widehat{P}$.
Notice that $G_T$ has a $\{1,1\}$-edge in $A_2$.
After putting negative loops at $u_{p/2}$ and $u_p$, we cannot locate a positive loop at $u_1$.
\end{proof}

\begin{lemma}
$p=2$.
\end{lemma}

\begin{proof}
Assume not.
Suppose $p/2$ is odd.
Then $\sigma$ has at least two orbits.
Thus the edges of $A_1$ form at least two essential cycles on $\widehat{P}$, where
$u_{p/2}$ and $u_p$ lie on distinct orbits.
Notice that $A_2$ contains a $\{p/2,p\}$-edge.
Since $u_{p/2}$ and $u_p$ are incident to negative loops, the edges of $A_1$ form just two cycles,
which are separating on $\widehat{P}$.
Furthermore, $u_1$ and $u_{p-1}$ lie one the same cycle.
Although there is a $\{1,p-1\}$-loop among positive loops at $v_1$,
we cannot locate it in $G_P$.

Suppose $p/2$ is even.
Then $\sigma$ has a single orbit.
Thus the edges of $A_1$ form an essential cycle $C$ on $\widehat{P}$.
Notice that $G_T$ contains a $\{1,p/2+1\}$-edge $e$ in $A_2$ and a $\{1,p-1\}$-loop $f$ at $v_1$.
After putting the negative loops at $u_{p/2}$ and $u_p$, 
we cannot locate $e$ (resp.\ $f$) in $G_P$ when $C$ is non-separating (resp.\ separating) on $\widehat{P}$.
\end{proof}

Finally, we eliminate the case $p=2$.
We denote by $\sigma$ the associated permutation to $A_1$.

\begin{lemma}\label{lem:2121}
If $G_T\cong G(2,2,1,2,1)$, then $\partial M$ is a single torus.
\end{lemma}

\begin{proof}
If $\sigma$ is the identity, then 
each vertex of $G_P$ is incident to two positive loops and two negative loops.
Hence these positive loops are separating on $\widehat{P}$.
Also, $G_T$ has two negative $\{1,2\}$-edges.
There are two possibilities for the arrangement of these two edges in $G_P$.
But both contradict Lemma \ref{lem:jumping} by looking the endpoint
of the edge of $A_2$ at $u_1$.

Thus $\sigma=(12)$.
Each vertex of $G_P$ is incident to one positive loop and two negative loops, and there are
$4$ positive edges between $u_1$ and $u_2$.
Then $G_P^+$ is contained in an annulus, whose core is separating on $\widehat{P}$.
By Lemma \ref{lem:common}, the $4$ positive edges between $u_1$ and $u_2$ are divided into two edge classes.
Then the jumping number is two, and $G_P$ is uniquely determined.
Let $M(\alpha)=N(\widehat{P})\cup W$.
Then $T\cap W$ consists of four bigons and four $3$-gons.
Let $D_1$ be the bigon contained in the parallelism between two loops at $v_1$, and $D_2$ 
a bigon between $v_1$ and $v_2$.
Also, let $D_3$ be any 3-gon.
Then $X=N(\partial W\cup (V_\alpha\cap W)\cup D_1\cup D_2\cup D_3)$ has $\partial W$ and a $2$-sphere as its boundary.
Since $M(\alpha)$ is irreducible, this implies $M(\alpha)$ is closed.
Hence $\partial M$ is a single torus.
\end{proof}

\begin{lemma}
$G_T\cong G(2,2,1,1,2)$ is impossible.
\end{lemma}

\begin{proof}
By the same argument as in the proof of Lemma \ref{lem:2121}, $\sigma=(12)$.
Again, the $4$ edges between $u_1$ and $u_2$ are divided into two edge classes.
In fact, they form two $S$-cycles, whose faces lie on the same side $Y$ of $\widehat{T}$.
By examining the edge correspondence, we see that the jumping number is two.
But we cannot draw two loops of the faces of those $S$-cycles on a genus two surface obtained from
$\widehat{T}$ by tubing along $V_\beta\cap Y$, simultaneously.
\end{proof}

\subsection{$q_1=(p+1)/2$}

\begin{lemma}
The case that $q_1=(p+1)/2$ is impossible.
\end{lemma}

\begin{proof}
Since $q_2+q_3+q_4+q_5=4p-1$, $G_T\cong G((p+1)/2,p,p,p,p-1)$.
Then two families $A_2$ and $A_3$ associate to the same permutation $\sigma$.
By Lemma \ref{lem:involution}, $\sigma^2$ must be the identity,
but this is impossible, because $p$ is odd.
Thus $\sigma$ is the identity.
Also, if $\tau$ is the associated permutation to $A_4$, then $\tau(i)\equiv i-1 \pmod{p}$.
Hence the edges in $A_4$ form an essential orientation-preserving cycle on $\widehat{P}$.
Since any vertex is incident to a positive loop, corresponding to the edges of $A_2$,
$G_P$ would contain a trivial loop.
\end{proof}

%
%
\section{Reduced graphs}\label{sec:reduced}

In this section, we prepare some results about the reduced graph of a graph $G$ (or its subgraph) on a Klein bottle $\widehat{P}$, which will be needed in the last section.
We need only the assumption that $G$ has no trivial loops and that the edges of $G$ are divided into positive edges and negative edges.

Let $\Lambda$ be a component of $\overline{G}^+$.
If there is a disk $D$ in $\widehat{P}$ such that $\mathrm{Int}\,D$ contains $\Lambda$, then we say that $\Lambda$ has a \textit{disk support}.
Also, if there is an annulus $A$ in $\widehat{P}$ such that $\mathrm{Int}\,A$ contains $\Lambda$ and $\Lambda$ does not have a disk support,
then we say that $\Lambda$ has an \textit{annulus support}.

Now, suppose that $\Lambda$ has a support $E$, where $E$ is a disk or an annulus.
A vertex $x$ of $\Lambda$ is called an \textit{outer vertex\/} if there is an arc $\xi$ connecting $x$ to $\partial E$
whose interior is disjoint from $\Lambda$.
Define an \textit{outer edge\/} similarly.
Then $\partial \Lambda$ denotes the subgraph of $\Lambda$ consisting of all outer vertices and all outer edges of $\Lambda$.
A vertex $x$ of $\Lambda$ is called a \textit{cut vertex\/} if $\Lambda-x$ has more components than $\Lambda$.

Suppose that $\Lambda$ has an annulus support $A$.
A vertex $x$ of $\Lambda$ is a \textit{pinched vertex\/} if
there is a spanning arc of $A$ which meets $\Lambda$ in only $x$.
An edge $e$ of $\Lambda$ is a \textit{pinched edge\/} if
there is a spanning arc of $A$ which meets $\Lambda$ in only one point on $e$.
Clearly, both endpoints of a pinched edge are pinched vertices.

We say that $\Lambda$ is an \textit{extremal component\/} of $\overline{G}^+$ if $\Lambda$ has a support which is disjoint from
the other components of $\overline{G}^+$.

\begin{lemma}\label{lem:support}
$\overline{G}^+$ has an extremal component with a disk support or an annulus support.
\end{lemma}

\begin{proof}
Let $\Lambda$ be a component of $\overline{G}^+$.
Choose a spanning tree $H$ of $\Lambda$, and
contract $H$ into one point.
Then we get a bouquet $\Lambda'$ in $\widehat{P}$.
Note that any loop in $\Lambda'$ is orientation-preserving.
If all loops in $\Lambda'$ are inessential in $\widehat{P}$, then
$\Lambda'$ has a disk support.
There are two isotopy classes of orientation-preserving essential loops in $\widehat{P}$.
But these two classes cannot exist simultaneously.
Therefore, if some loop in $\Lambda'$ is essential,
then $\Lambda'$ has an annulus support, and so does $\Lambda$.

If $\overline{G}^+$ has a component with a disk support,
then there exists an extremal component with a disk support.
Otherwise, any component of $\overline{G}^+$ has an annulus support, and hence
any component is extremal.
\end{proof}

Let $x$ be a vertex of $G$.
Then $x$ is called an \textit{interior vertex\/} if there is no negative edge incident to $x$ in $G$.
Since $\overline{G}$ and $\overline{G}^+$  have the same vertex set as $G$,
we may call a vertex of $\overline{G}$ or $\overline{G}^+$ an interior vertex when it is an interior vertex of $G$.
In particular, if $x$ is in an extremal component of $\overline{G}^+$ with a disk or an annulus support, and
it is not an outer vertex, then $x$ is an interior vertex.

A vertex $x$ is said to be \textit{good\/} if all positive edge endpoints around $x$ are successive in $G$.
Thus an interior vertex is good.
When $x$ is a vertex of an extremal component $\Lambda$ of $\overline{G}^+$ with a disk or an annulus support,
$x$ is good if
\begin{itemize}
\item[(i)] $x$ is not a cut vertex of $\Lambda$ if $\Lambda$ has a disk support; or
\item[(ii)] $x$ is neither a cut vertex nor a pinched vertex of $\Lambda$ if $\Lambda$ has an annulus support.
\end{itemize}

\begin{proposition}\label{prop:int}
If each interior vertex of $\overline{G}$ has degree at least $6$, then
$\overline{G}^+$ has either a good vertex of degree at most $4$, or a vertex of degree at most $2$.
\end{proposition}

\begin{proof}
See \cite[Proposition 3.4]{T} (and its proof).
If an extremal component of $\overline{G}^+$ is a single vertex or a cycle, then
we have the second conclusion.
\end{proof}

If $G$ has no interior vertex, then we have a stronger conclusion.

\begin{figure}[htb]
\begin{center}
\includegraphics*[scale=0.27]{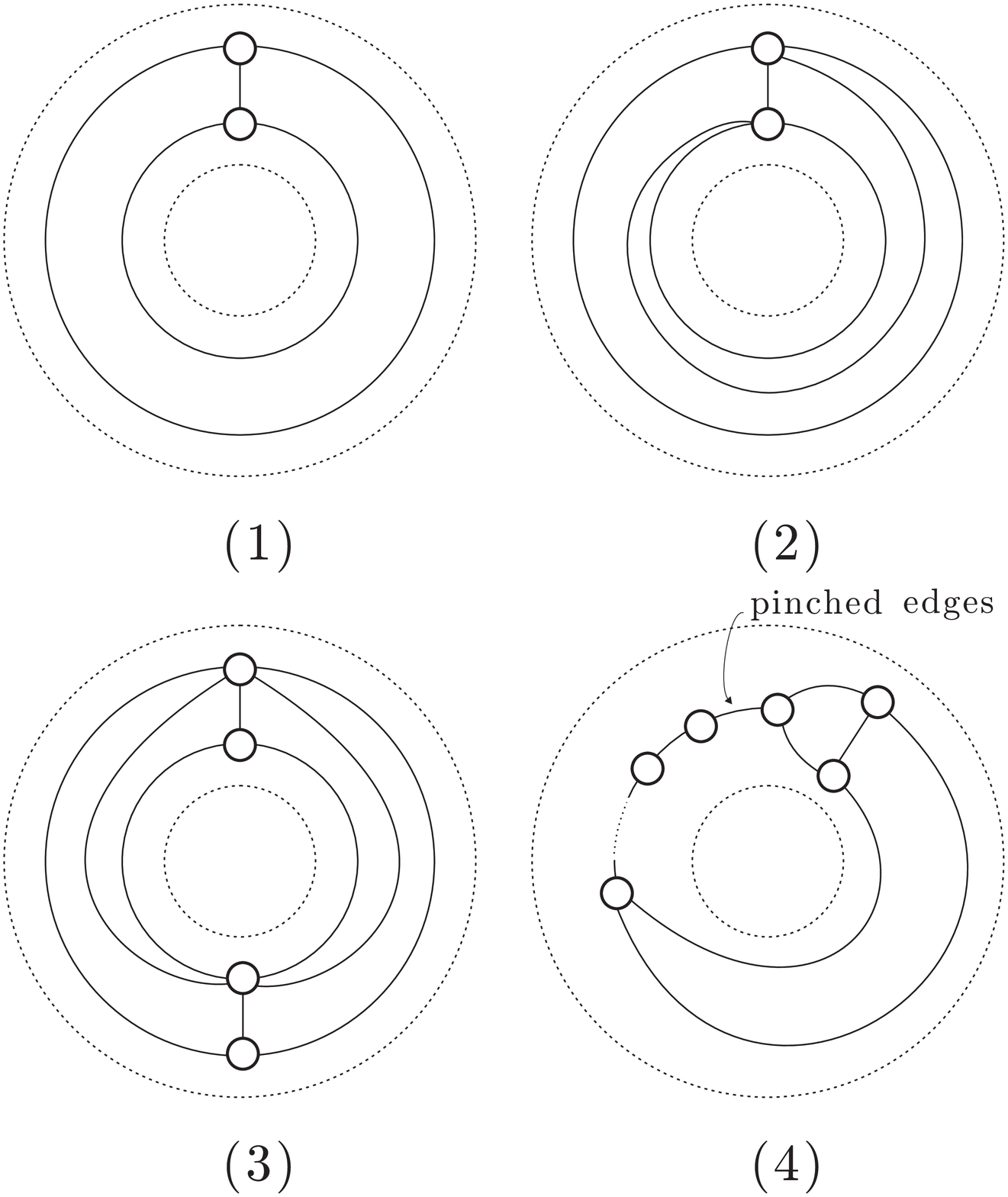}
\caption{}\label{fig:a-supp}
\end{center}
\end{figure}

\begin{lemma}\label{lem:noint-a}
Suppose that $G$ has no interior vertex.
Let $\Lambda$ be an extremal component of $\overline{G}^+$.
If $\Lambda$ has an annulus support and $\Lambda$ is not a cycle, then either
\begin{itemize}
\item[\rm(1)] $\Lambda$ has two non-pinched good vertices of degree at most $4$ on the same side of $\Lambda$\textup{;}
\item[\rm(2)] $\Lambda$ has a non-pinched vertex of degree at most two\textup{;} or
\item[\rm(3)] $\Lambda$ is as shown in Figure \textup{\ref{fig:a-supp}(1)}, \textup{(2)}, \textup{(3)} or \textup{(4)} with possibly no pinched edge.
\end{itemize}
\end{lemma}

\nocolon

\begin{proof}
Let $V$ be the number of vertices of $\Lambda$.

(1)\qua First, consider the case where $\Lambda$ has no cut vertex.

If $\Lambda$ has no pinched vertex, then $\partial \Lambda$ consists of two cycles.
Note that any vertex lies on $\partial \Lambda$, because $G$ has no interior vertex.
If some edge, not in $\partial \Lambda$, connects two vertices on the same cycle of $\partial \Lambda$,
then the cycle contains a vertex of degree two, giving the conclusion (2).
Hence we can assume that any edge of $\Lambda$, not in $\partial \Lambda$, connects two vertices on distinct sides.

If $V=2$, then $\Lambda$ is as shown in  Figure \ref{fig:a-supp}\textup{(1)} or \textup{(2)}.
If $V=3$, then it is easy to see that (1) or (2) holds.
Let $V=4$.
If one cycle of $\partial\Lambda$ contains one vertex, then we have (1) or (2).
Hence we may assume that each cycle of $\partial \Lambda$ contains two vertices.
Clearly, any vertex has degree at most $5$.
If all vertices have degree at most $4$, then we have (1).
So, assume that some vertex has degree $5$.
Then we have (2), or $\Lambda$ is Figure \ref{fig:a-supp}\textup{(3)}.
Hereafter we assume $V\ge 5$.

If there are more than two vertices of degree at most $4$, then we have (1).
Hence we assume that all vertices but at most two vertices $x, y$ have degree at least $5$.
Take a double of $\Lambda$ along $\partial \Lambda$, and let $E$ and $F$ be the number of edges, faces, respectively,
as a graph on a torus.
Then $V-E+F=0$ and $3F\le 2E$, giving $E\le 3V$.
For a vertex of $\Lambda$ with degree at least $5$, it has degree at least $8$ in the double.
Hence $\mathrm{deg}(x)+\mathrm{deg}(y)+8(V-2)\le 2E\le 6V$, where $\mathrm{deg}(-)$ denotes degree in the double.
Then $\mathrm{deg}(x)+\mathrm{deg}(y)\le 16-2V\le 6$.
If either $x$ or $y$ has degree two in $\Lambda$, then (2) holds.
But, if not, $\mathrm{deg}(x)\ge 4$ and $\mathrm{deg}(y)\ge 4$, a contradiction.

Next, suppose that $\Lambda$ has a pinched vertex.
If necessary, contract all pinched edges, and denote the resulting graph by $\Lambda'$.
If $\Lambda'$ contains more than one pinched vertices, then consider a part $H$ between two consecutive pinched vertices.
If $\Lambda'$ contains only one pinched vertex $x$, then
take a spanning arc $\xi$ of the annulus support of $\Lambda'$ with $\xi\cap \Lambda'=x$, and
split along $\xi$ to obtain $H$.
In any case, $H$ has more than two vertices, and has a disk support.
Let $x_1$ and $x_2$ denote the two vertices coming from pinched vertices of $\Lambda'$.

If $H=\partial H$, then any vertex, except $x_1$ and $x_2$, gives the conclusion (2).
Otherwise, $H$ has an edge not on $\partial H$.
If there is an edge incident to $x_i$ not on $\partial H$, then $H$ contains a good vertex of degree two by considering
an outermost edge.
Thus we may assume that $\mathrm{deg}(x_1)=\mathrm{deg}(x_2)=2$.
Let $V'$, $E'$, $F'$ be the numbers of vertices, edges and faces of $H$ as a graph in a disk.
Then $V'\ge 4$, and $V'-E'+F'=1$, $3F'+V'\le 2E'$, giving $E'\le 2V'-3$.

\begin{claim}\label{cl:h}
If $V'>4$, then $H$ has either a good vertex of degree two, or at least $3$ good vertices of degree at most $4$, except $x_1$ and $x_2$.
\end{claim}

\begin{proof}[Proof of Claim \ref{cl:h}]
Assume that any vertex, except $x_1$, $x_2$, $y$, $z$ has degree at least $5$.
Then $\mathrm{deg}(x_1)+\mathrm{deg}(x_2)+\mathrm{deg}(y)+\mathrm{deg}(z)+5(V'-4)\le 2E'$.
So, $\mathrm{deg}(y)+\mathrm{deg}(z)\le 10-V'\le 5$.
Thus $y$ or $z$ has degree two.
\end{proof}

Thus if $V'>4$, then we have the conclusion (1) or (2).
When $V'=4$, the other two vertices of $H$ than $x_1$ and $x_2$ are connected with a single edge, and so have degree $3$.
If $\Lambda'$ contains more than one pinched vertices, then there are at least two parts such as $H$.
Then we have the conclusion (1).
Otherwise, $\Lambda$ must be the form as in Figure \ref{fig:a-supp}(4).

(2)\qua Consider the case where $\Lambda$ has a cut vertex.

If some block has a disk support, then we can see that (2) holds by \cite[Lemma 3.2]{W}.
Thus we can assume that any block has an annulus support.
Then either $\partial \Lambda$ consists of two cycles, or $\Lambda$ has a single pinched vertex.
Hence the first or the second part of the previous case gives the result, respectively.
\end{proof}

\begin{figure}[htb]
\begin{center}
\includegraphics*[scale=0.30]{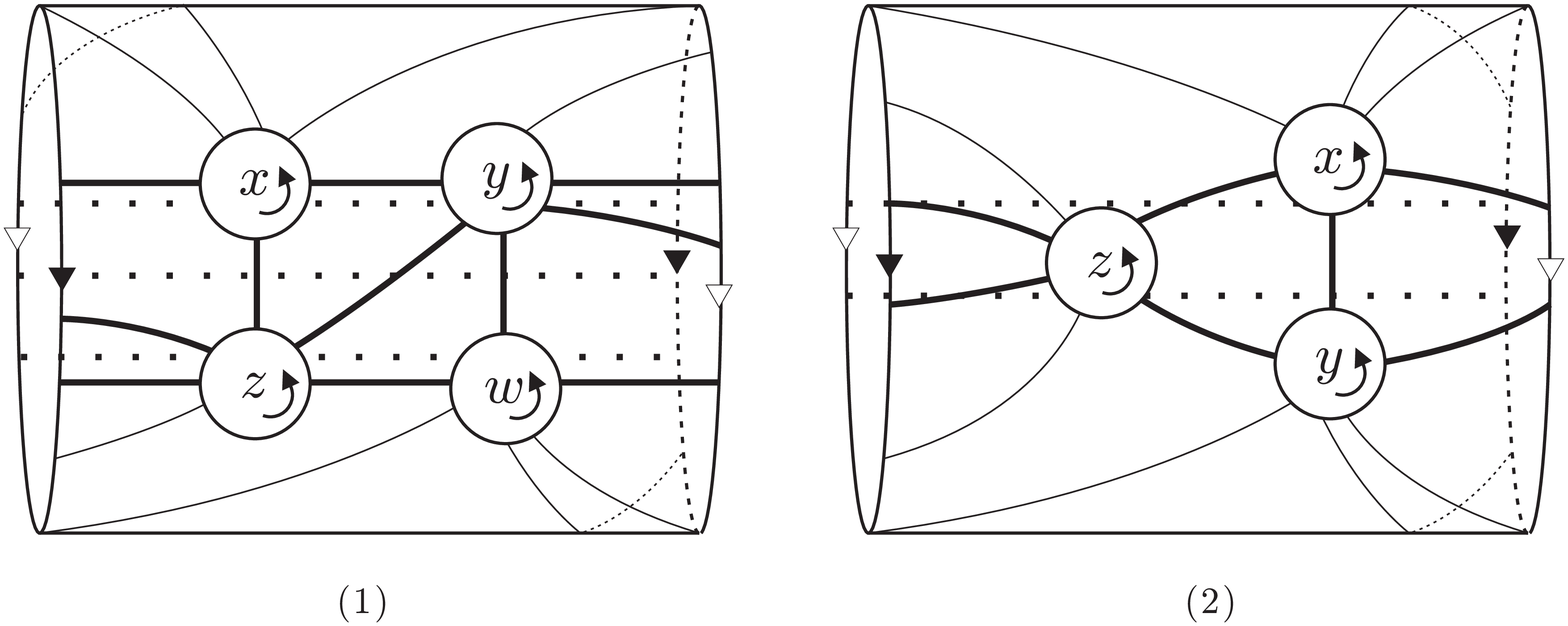}
\caption{}\label{fig:special3}
\end{center}
\end{figure}

\begin{proposition}\label{prop:klein}
Suppose that $G$ lies on a Klein bottle and has $p\ge 3$ vertices.
If $G$ has no interior vertex, then either
\begin{itemize}
\item[\rm(1)] $\overline{G}^+$ has a good vertex of degree $4$, which is not incident to a negative loop in $G$\textup{;}
\item[\rm(2)] $\overline{G}^+$ has a vertex of degree at most $3$, which is not incident to a negative loop in $G$\textup{;}
\item[\rm(3)] $\overline{G}^+$ has a vertex of degree at most $2$, which is incident to a single negative loop in $\overline{G}$\textup{;} or
\item[\rm(4)] $\overline{G}$ is either of the graphs
shown in Figure \textup{\ref{fig:special3}}, where the end circles of the cylinder are identified suitably to form a Klein bottle,
and the thicker edges are positive.
\end{itemize}
\end{proposition}

\nocolon

\begin{proof}
Let $\Lambda$ be an extremal component of $\overline{G}^+$.
Assume that $\Lambda$ has a disk support.
If $\Lambda$ is not a single vertex, then it has two good vertices of degree at most two \cite[Lemma 3.2]{W}.
If either of them is not incident to a negative loop, then (2) holds.
Otherwise, either we have (3), or
both are incident to more than one negative loops in $\overline{G}$.
Then there is another extremal component $\Lambda'$ with a disk support.
Notice that at most two vertices of $G$ can be incident to a negative loop.
Hence any vertex of $\Lambda'$ of degree at most two is not incident to a negative loop, which gives (2) again.
If $\Lambda$ is a single vertex, then either we have (2) or (3), or $\Lambda$ has more than one negative loops in $\overline{G}$.
But the latter implies the existence of another extremal component with a disk support, which gives (2) as above.

Thus we assume that $\overline{G}^+$ has no component with a disk support.
Hence any component has an annulus support by Lemma \ref{lem:support}, and anyone is extremal.
First, if all components of $\overline{G}^+$ are cycles,
we can choose a vertex of degree two, which is not incident to a negative loop, since $p\ge 3$.
So we have (2).
Thus we may choose so that $\Lambda$ is not cycle.
Then one of (1), (2), (3) in Lemma \ref{lem:noint-a} holds for $\Lambda$.

If (1) of Lemma \ref{lem:noint-a} happens, then either of such two vertices is not incident to a negative loop.
Therefore we have the conclusion (1) or (2).
If (2) of Lemma \ref{lem:noint-a} happens, then a non-pinched vertex of degree at most two in $\Lambda$
satisfies (2) or (3).
Finally, assume that $\Lambda$ satisfies (3) of Lemma \ref{lem:noint-a}. 
If $\Lambda$ is as Figure \ref{fig:a-supp}\textup{(1)} or \textup{(2)}, then
it has a good vertex of degree at most $4$, which is not incident to a negative loop, because $\overline{G}^+$ has other component.
This gives the conclusion (1) or (2).

If $\Lambda$ is as Figure \ref{fig:a-supp}\textup{(3)} and $\overline{G}^+$ has other component, then
$\Lambda$ has a good vertex of degree $3$, which is not incident to a negative loop.
This is (2).
Hence we assume $\overline{G}^+=\Lambda$ has the form Figure \ref{fig:a-supp}\textup{(3)}.
Let $A$ be its annulus support.
If the core of $A$ is non-separating on the Klein bottle, we have (2) again.
Hence the core of $A$ is assumed to be separating.
Let $x$ and $w$ be the vertices of degree $3$ in $\Lambda$.
If either of them is not incident to a negative loop, then (2) holds.
Hence we assume that both are incident to a negative loop.
Since there is no interior vertex, $\overline{G}$ must be the graph of Figure \ref{fig:special3}(1).

If $\Lambda$ is as Figure \ref{fig:a-supp}\textup{(4)} and $\overline{G}^+$ has other component, then
$\Lambda$ has a good vertex of degree $3$, which is not incident to a negative loop.
This is (2) again.
So, suppose that $\overline{G}^+=\Lambda$ has the form Figure \ref{fig:a-supp}\textup{(4)}.
Let $x$ and $y$ be the non-pinched vertices of $\Lambda$ of degree $3$.
If either of them is not incident to a negative loop, then (2) holds.
Assume that both are incident to a negative loop.
If $p>3$, then any pinched vertex satisfies (2).
If $p=3$, then the unique pinched vertex satisfies (1) or $\overline{G}$ is the graph of Figure \ref{fig:special3}(2).
\end{proof}

\begin{figure}[htb]
\begin{center}
\includegraphics*[scale=0.33]{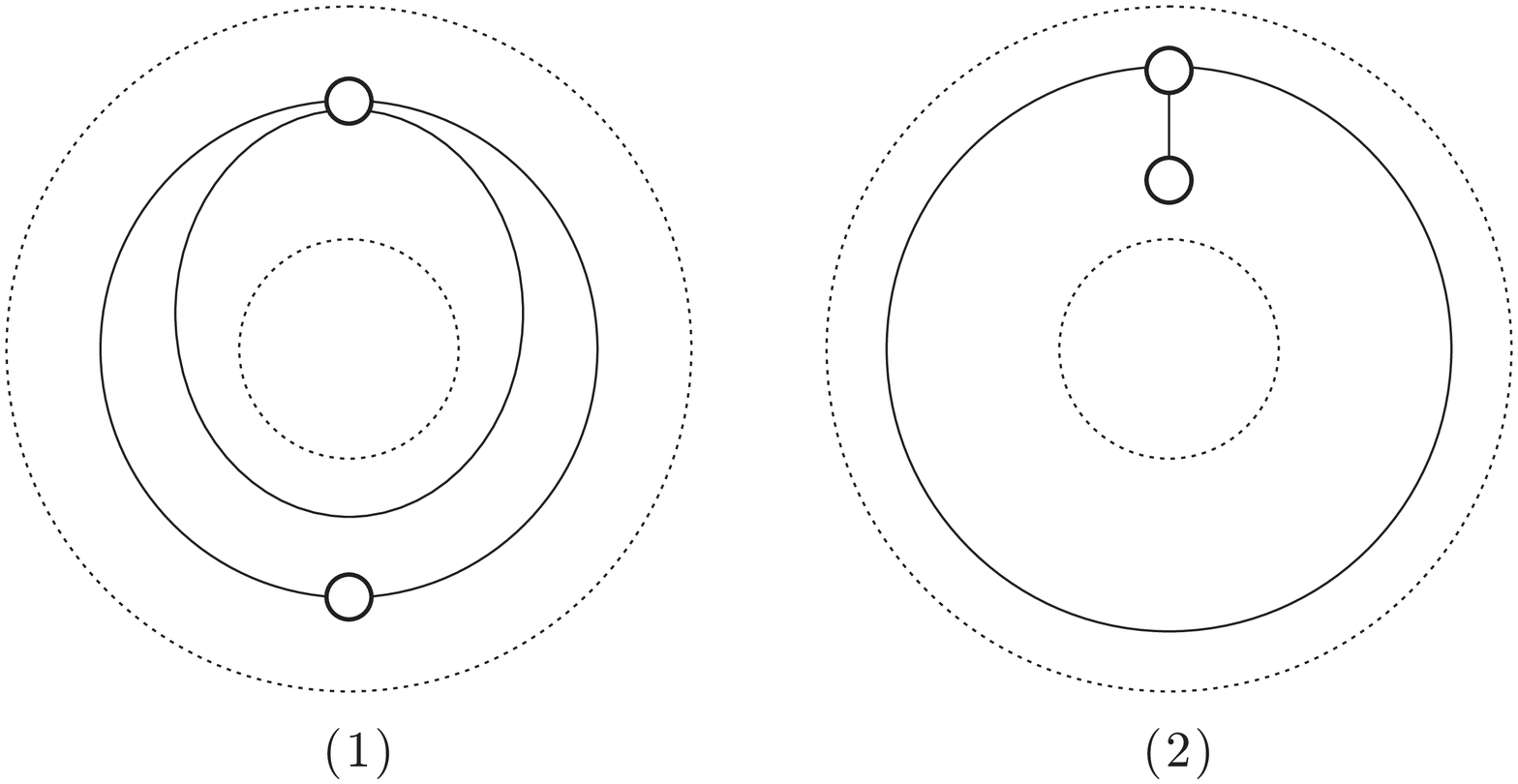}
\caption{}\label{fig:2vertex}
\end{center}
\end{figure}

\begin{lemma}\label{lem:2vertex}
If $G$ has only two vertices, then $\overline{G}^+$ is one of the following\textup{:}
\begin{itemize}
\item[\rm(1)] a single edge\textup{;}
\item[\rm(2)] a cycle of length two\textup{;}
\item[\rm(3)] the graph as shown in Figure \textup{\ref{fig:a-supp}(1)}, \textup{(2)} or Figure \textup{\ref{fig:2vertex}(1)} or \textup{(2)}\textup{;}
\item[\rm(4)] two isolated vertices\textup{;}
\item[\rm(5)] two loops\textup{;} or
\item[\rm(6)] an isolated vertex and a loop.
\end{itemize}
\end{lemma}

\nocolon

\begin{proof}
Let $u_1$ and $u_2$ be the vertices of $G$.
Notice that the number of loops in $G$ at $u_1$ is equal to that of loops at $u_2$.

Assume that $\overline{G}^+$ is connected.
If $\overline{G}^+$ has a disk support, then it is a single edge.
For, if there is a loop at $u_1$ say, then there would be a trivial loop at $u_2$.
If $\overline{G}^+$ has an annulus support, then it is easy to see that
(2) or (3) holds.

Next assume that $\overline{G}^+$ is not connected.
Let $H_1$, $H_2$ be the components of $\overline{G}^+$ containing $u_1$ and $u_2$, respectively.
If $H_1$ has a disk support, then $H_1=u_1$.
Also if $H_1$ has an annulus support, then $H_1$ is a positive loop.
Thus either (4), (5) or (6) holds.
\end{proof}

\section{Klein bottle; generic case}\label{sec:generic-kb}

Finally, we consider the case that $p\ge 2$ and $t\ge 3$.

\begin{lemma}\label{lem:t4key}
Assume $t=4$ and that $G_T$ contains
$m$ $x$-edges connecting $v_1$ and $v_2$ and $n$ $x$-edges connecting $v_3$ and $v_4$ for some label $x$.
Then
\begin{itemize}
\item[\rm(1)] if $m, n \ge 4$, then $m=n=4$;
\item[\rm(2)] if all the $x$-edges are level and $m,n\ge 2$, then $m=n=2$.
\end{itemize}
\end{lemma}

\begin{proof}
Let $G(1,2)$ be the subgraph of $G_T$
consisting of $v_1,v_2$ and $m$ $x$-edges between them.
Define $G(3,4)$ similarly.

(1)\qua
If $G(1,2)$, say, has a disk support disjoint from $G(3,4)$ on $\widehat{T}$,
then $G(1,2)$ contains three mutually parallel $x$-edges.
But this means that $G_T$ has $p+1$ parallel edges, which is impossible by Lemma \ref{lem:parallel-max3}.
Hence neither $G(1,2)$ nor $G(3,4)$ has a disk support on $\widehat{T}$.
Then both have annulus supports.
If $m>4$, then $G(1,2)$ would have three mutually parallel edges again, and so $m=4$.
Similarly, we have $n=4$.

(2)\qua Since two level $x$-edges cannot be parallel,
the result follows from a similar argument to (1).
\end{proof}

\begin{lemma}\label{lem:successive3}
If $t=4$, then
there are no three successive \textup{(}distinct\textup{)}
positive edges of weight $4$ incident to a vertex
of $\overline{G}_P$.
\end{lemma}

\begin{proof}
Let $e_1,e_2,e_3$ be successive positive edges incident to a vertex $x$ of $\overline{G}_P$.
If each $e_i$ has weight $4$, then we may assume that each family of mutually parallel
positive edges corresponding to $e_i$ contains a $\{1,2\}$ $S$-cycle and a $\{3,4\}$ $S$-cycle
by Lemma \ref{lem:parallel-max-kb}.
Then this contradicts Lemma \ref{lem:t4key}(1).
\end{proof}

\begin{lemma}\label{lem:key}
Let $u_i$ be a vertex of $G_P$.
Suppose that $u_i$ is incident to $k$ non-loop negative edges  and $n$ negative loops in $G_P$
and that $u_i$ has degree $m$ in $\overline{G}_P^+$.
Then, 
\begin{itemize}
\item[\rm(1)] $k\le $t\textup{;}
\item[\rm(2)] $k+2n\ge (10-m)t/2-2m$ when $t\equiv 0\pmod{4}$, and 
$k+2n\ge (10-m)t/2-m$, otherwise\textup{;} and
\item[\rm(3)] If $\widehat{T}$ is non-separating, then $k+2n\ge (10-m)t/2$.
\end{itemize}
\end{lemma}

\begin{proof}
In $G_T$, there are $n$ positive level $i$-edges and $k$ positive non-level $i$-edges.
By Lemma \ref{lem:parallel-max3}, no two of positive $i$-edges are parallel.
Also, a negative level $i$-edge cannot be parallel to a negative $i$-edge, and a negative non-level $i$-edge
can be parallel to at most one negative $i$-edge.
Thus $\overline{G}_T$ has at least $n+k+(5t-(2n+k))/2=(5t+k)/2$ edges.
Since $\overline{G}_T$ has at most $3t$ edges, we have $(5t+k)/2\le 3t$, giving (1).
Also, $(k+2n)+m(t/2+2)\ge 5t$  and $(k+2n)+m(t/2+1)\ge 5t$ give (2).
If $\widehat{T}$ is non-separating, then $(k+2n)+m\cdot t/2\ge 5t$, giving (3).
\end{proof}

\begin{proposition}\label{prop:key}
Suppose that $G_P$ contains a positive edge.
Let $u_i$ be a vertex of $G_P$, which has degree $m$ in $\overline{G}_P^+$.
Assume that $u_i$ is incident to at most one negative loop in $\overline{G}_P$, and
let $n\,(\ge 0)$ denote the weight of the negative loop.
If $n=0$, then $m\ge 4$, and otherwise, $m\ge 2$.
When the equality holds, $\widehat{T}$ is separating and $t=4$.
\end{proposition}

\begin{proof}
By Lemma \ref{lem:parallel-max-kb}, $n\le t$.
Then the conclusion immediately follows from Lemma \ref{lem:key}.
\end{proof}


\subsection{Case 1: $p\ge 3$}

\begin{proposition}\label{prop:six}
Each vertex of $\overline{G}_T$ has degree $6$.
\end{proposition}

\begin{proof}
Let $v_i$ be a vertex of $\overline{G}_T$.
If $\mathrm{deg}(v_i)\ge 6$, then an easy Euler characteristic calculation shows
$\mathrm{deg}(v_i)=6$.
Since each edge has weight at most $p$ by Lemma \ref{lem:parallel-max3},
$\mathrm{deg}(v_i)\ge 5$.
Hence suppose $\mathrm{deg}(v_i)=5$ for contradiction.
Then $v_i$ is incident to $5$ negative edges with weight $p$, because $p/2+1<p$.
By the parity rule, all $i$-edges of $G_P$ are positive.
In particular, any vertex of $G_P$ is incident to five $i$-edges, none of which are parallel by Lemma \ref{lem:parallel-max-kb}.
Thus any vertex of $\overline{G}_P^+$ has degree at least $5$.
By Proposition \ref{prop:int}, 
$\overline{G}_P$ has an interior vertex $u$ of degree at most $5$, hence just $5$.
Then $5(t/2+2)\ge 5t$ gives $t\le 4$.
But, if $t=3$, then $G_P$ cannot contain a pair of parallel positive edges.
So $t=4$.
Thus $u$ is incident to $5$ positive edges of weight $4$.
In $G_P$, each edge corresponds to a family of $4$ parallel positive edges.
By Lemma \ref{lem:parallel-max-kb}(1), we may assume that each family contains a $\{1,2\}$ $S$-cycle and a $\{3,4\}$ $S$-cycle.
Then we can see that $u$ is not incident to a loop.
But this contradicts Lemma \ref{lem:successive3}.
\end{proof}

\begin{lemma}
$G_P$ contains a positive edge.
\end{lemma}

\begin{proof}
If not, we can choose a vertex which is not incident to a negative loop, because
$p\ge 3$ and at most two vertices can be incident to a negative loop.
But this contradicts Lemma \ref{lem:key}(1).
\end{proof}

We now consider two cases, according to the existence of an interior vertex in $G_P$.

\subsubsection{Case: $G_P$ has no interior vertex}

\begin{lemma}\label{lem:degree4}
$\overline{G}_P^+$ cannot have a good vertex of degree $4$, which is not incident to a negative loop in $G_P$.
\end{lemma}

\begin{proof}
Let $u_i$ be such a vertex, and let $k$ be the number of negative edge endpoints at $u_i$.
By Proposition \ref{prop:key} and Lemma \ref{lem:key}, $t=k=4$.
Hence $u_i$ is incident to $4$ families of $4$ mutually parallel positive edges, successively, because
$u_i$ is good.
By examining the labels, we see that there is no positive loop at $u_i$.
But this contradicts Lemma \ref{lem:successive3}.
\end{proof}

\begin{lemma}\label{lem:degree3}
$\overline{G}_P^+$ cannot have a vertex of degree at most $3$, which is not incident to a negative loop in $G_P$.
\end{lemma}

\begin{proof}
This is an immediate consequence of Proposition \ref{prop:key}.
\end{proof}

\begin{lemma}\label{lem:degree2}
$\overline{G}_P^+$ cannot have a vertex of degree at most two, which is incident to a single negative loop in $\overline{G}_P$. 
\end{lemma}

\begin{proof}
Let $u_i$ be such a vertex.
Let $n$ be the number of negative loops and $k$ be the number of negative non-loop edges at $u_i$ in $G_P$.
By Proposition \ref{prop:key}, $\widehat{T}$ is separating and $t=4$.
Then $u_i$ has at most $8$ positive edge endpoints, and so at least $12$ negative edge endpoints.
Notice that $G_T$ has $n$ positive level $i$-edges and at least $12-2n$ positive non-level $i$-edges, and that
no two of them are parallel.
Hence $\overline{G}_T$ has at least $n+(12-2n)=12-n\ge 8$ positive $i$-edges, and then at most $4$ negative edges.

Let $\Lambda$ be an extremal component of $\overline{G}_T^+$.
If $\Lambda$ has a disk support, then we see that it contains a vertex of degree at most one.
Then such a vertex is incident to at least $5$ negative edges in $\overline{G}_T$.
But this is impossible, because $\overline{G}_T$ has at most $4$ negative edges.
Thus $\overline{G}_T^+$ has no component with a disk support, and then each component has an annulus support.

Let $G(1,3)$ be the (possibly, disconnected) subgraph of $\overline{G}_T^+$ spanned by $v_1$ and $v_3$, and
define $G(2,4)$ similarly.
If $G(1,3)$ is disconnected, then it consists of two loops.
Then $G(2,4)$ contains at least $6$ edges.
But this is impossible by an Euler characteristic calculation.
This implies that both of $G(1,3)$ and $G(2,4)$ are connected.
Then we see that both of them contain $4$ edges and have the form as in Figure \ref{fig:a-supp}(2).
Each vertex of $\overline{G}_T$ has at most two negative edges.
Thus $4(p/2+1)+2p\ge 5p$ gives $p\le 4$.

If $p=3$, then each positive edge at $v_1$ has weight at most two, and so the total weight cannot be $5p=15$.
If $p=4$, then the $4$ positive edges at $v_1$ have weight $3$ and two negative edges have weight $4$.
But there is an $S$-cycle among loops at $v_1$, which is impossible by Lemma \ref{lem:kb-S}.
\end{proof}


\begin{lemma}\label{lem:worst1}
$\overline{G}_P$ is not the graph as shown in Figure \textup{\ref{fig:special3}(1)}.
\end{lemma}

\begin{proof}
Let $x$ and $y$ be the vertices as shown there.
If $k$ is the number of negative edges at $y$ in $G_P$, then $k\le t$ by Lemma \ref{lem:key}.
For $y$, Lemma \ref{lem:key}(2) gives that $\widehat{T}$ is separating and $t=4$.
Let $n$ be the number of negative loops at $x$.
By Lemma \ref{lem:successive3}, $x$ has at most $11$ positive edge endpoints, and then $k+2n\ge 9$.
Hence $n=3$ or $4$.

Assume $n=3$.
Notice that $k$ must be odd by the parity rule.
Hence $k=3$.
Thus the three positive edges at $x$ have weights $\{4,4,3\}$.
By examining the labels at $x$, we see that all vertices of $G_T$ have degree at least $6$ in $\overline{G}_T$
and some vertex has degree more than $6$ there.
(For example, see Figure \ref{fig:worst1}.
In this case, the degrees of $v_i$ in $\overline{G}_T$ are $7$, $7$, $6$ and $6$ for $i=1,2,3,4$, respectively.)
Clearly, this is impossible.

\nocolon\begin{figure}[htb]
\begin{center}
\includegraphics*[scale=0.45]{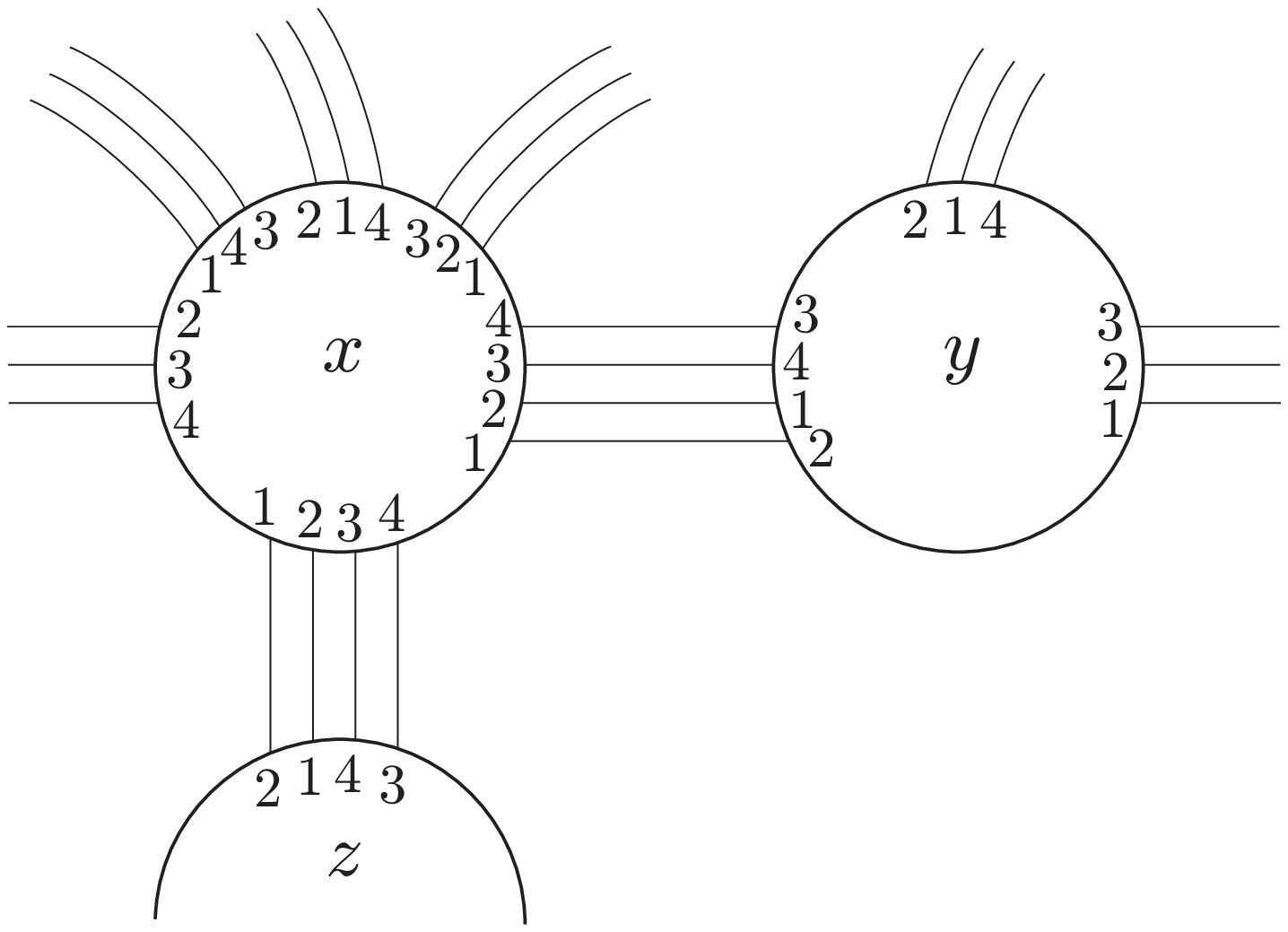}
\caption{}\label{fig:worst1}
\end{center}
\end{figure}

Assume $n=4$.   By the parity rule, $k$ must be even.
Then $k=2$ or $4$.

First assume $k=2$.
Then the three positive edges at $x$ have weights $\{4,4,2\}$ or $\{4,3,3\}$.
If one family between $x$ and $y$ has weight $4$, then
look at the two vertices, which do not appear at non-loop negative edge endpoints at $x$.
In $\overline{G}_T$, they cannot have degree $6$, which contradicts Proposition \ref{prop:six}.
Thus the only possibility is that the two edges between $x$ and $y$ have weight $3$, and the edge between $x$ and $w$ has weight $4$.
Then the two non-loop edges at $x$ are level by Lemma \ref{lem:t4key}(1).
But two vertices of $\overline{G}_T$ cannot have degree $6$ as above again.

Thus we have $k=4$.
Then the associated permutation to the family of $4$ negative loops at $x$ is the identity.
Hence any vertex of $G_T$ is incident to a loop.
If the associated permutation $\sigma$ to the family of $4$ non-loop negative edges at $x$ is also the identity,
then a level $x$-loop and a $\{x,y\}$-loop are incident to each vertex of $G_T$.
These two loops are parallel, which contradicts Lemma \ref{lem:parallel-max3}.
Thus $\sigma$ is the permutation $(13)(24)$.
We may assume that the labels in $G_P$ are as in Figure \ref{fig:worst1-2}, where
$a$, $b$ and $c$ denotes the number of edges in the families.
Hence $\overline{G}_T^+$ has two components, each of which has the form as in Figure \ref{fig:a-supp}(2).
Notice that $G_P$ contains a $\{1,2\}$-edge and a $\{3,4\}$-edge.
This implies that any negative edge of $G_T$ connects either $v_1$ and $v_2$, or $v_3$ and $v_4$.
By the parity rule, any positive edge of $G_P$ is a $\{1,2\}$-edge or a $\{3,4\}$-edge.

\nocolon\begin{figure}[htb]
\begin{center}
\includegraphics*[scale=0.45]{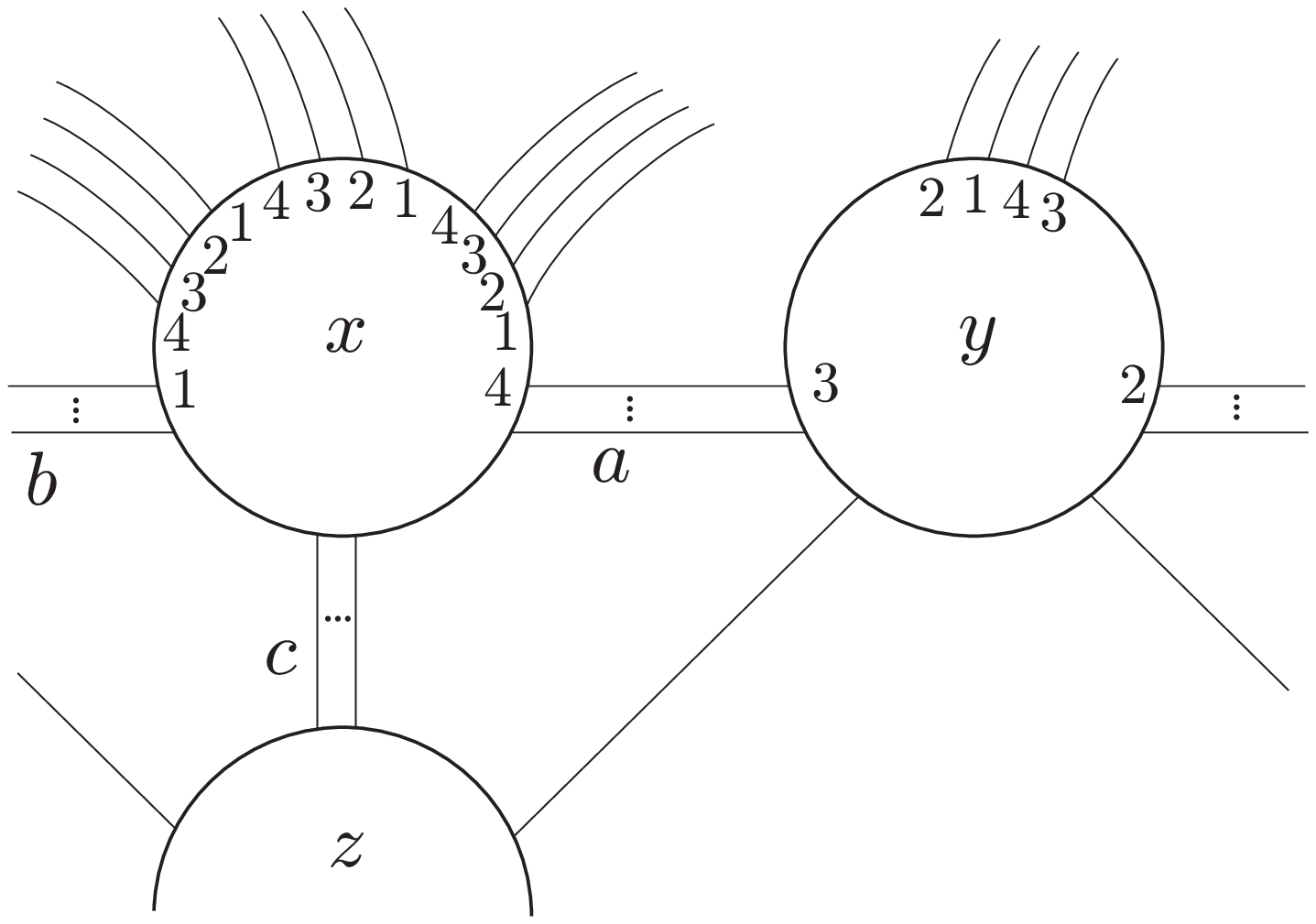}
\caption{}\label{fig:worst1-2}
\end{center}
\end{figure}

If $a=2$ or $4$, then there would be a $\{1,4\}$- or $\{2,3\}$-edge between $y$ and $z$, a contradiction.
Similarly, $b\ne 2$, $4$.
Since $a+b+c=8$ and none of them is zero, $(a,b,c)=(1,3,4)$, $(3,1,4)$ or $(3,3,2)$.
If $(a,b,c)=(1,3,4)$ or $(3,1,4)$, then the family between $x$ and $z$ contains an extended $S$-cycle, a contradiction.
Finally, if $(a,b,c)=(3,3,2)$, then $G_P$ contains $S$-cycles with label set $\{1,2\}$ and $\{3,4\}$ and
a Scharlemann cycle of length $3$ with label set $\{1,2\}$.
But, under the existence of a Scharlemann cycle with label set $\{3,4\}$,
all Scharlemann cycles with label set $\{1,2\}$ must have the same length \cite[Theorem 5.7]{GL4}.
\end{proof}

\begin{lemma}\label{lem:worst2}
$\overline{G}_P$ is not the graph as shown in Figure \textup{\ref{fig:special3}(2)}.
\end{lemma}

\begin{proof}
Let $x$ and $y$ be the vertices which are incident to a negative loop, and $z$ the other one.
By applying Lemma \ref{lem:key} and Proposition \ref{prop:key} to $z$, $t=k=4$, where $k$ is the number of negative edge at $z$ in $G_P$.
Then each positive edge at $z$ has weight $4$.
Since $z$ is not good, the weights of two negative edges at $z$ are $\{1,3\}$ or $\{2,2\}$. 
By examining the labels, the former contradicts Lemmas \ref{lem:parallel-max-kb}(1) and \ref{lem:t4key}(1), and
the latter contradicts Lemma \ref{lem:parallel-max-kb}(1).
\end{proof}

\begin{proposition}\label{prop:exist-int}
$G_P$ must have an interior vertex.
\end{proposition}

\begin{proof}
If not, (1), (2), (3) or (4) of Proposition \ref{prop:klein} holds.
But all of these are impossible by Lemmas \ref{lem:degree4},
\ref{lem:degree3}, \ref{lem:degree2}, \ref{lem:worst1} and \ref{lem:worst2}.
\end{proof}

\subsubsection{Case: $G_P$ has an interior vertex}

Let $u_i$ be an interior vertex of $G_P$. 
Thus all edges incident to $u_i$ are positive, and hence all $i$-edges in $G_T$ are negative.

\begin{lemma}
$t=4$.
\end{lemma}

\begin{proof}
The argument of the proof of \cite[Lemma 4.5]{T} works without any change.
Therefore $G_P$ has an $S$-cycle with label $j$ for any label $j$.
In particular, $\widehat{T}$ is separating and $t$ is even.
If $t\ge 6$, then $G_P$ would contain three $S$-cycles with disjoint label pairs as in the proof of \cite[Proposition 4.6]{T}.
\end{proof}

\begin{lemma}\label{lem:4negative}
Any vertex of $G_P$ has at most $4$ negative edge endpoints.
\end{lemma}

\begin{proof}
Let $l$ be the number of (positive) loops at $u_i$ in $G_P$.
Then $G_T$ has $n$ negative level $i$-edges and $20-2l$ negative non-level $i$-edges.
Among these $i$-edges, none is parallel to a level one, and at most two non-level ones can be parallel.
Thus $\overline{G}_T$ has at least $l+(20-2l)/2=10$ negative edges, and hence at most two positive edges.

If a vertex $u_j$ has more than $4$ negative edge endpoints, then $u_j$ is incident to at least $3$ negative edges
(possibly, loops).
This means that $G_T$ has at least $3$ positive $j$-edges.
Thus $G_T$ contains a pair of parallel positive $j$-edges, and so an $S$-cycle or a generalized $S$-cycle.
This is impossible by Lemma \ref{lem:kb-S}.
\end{proof}

\begin{proposition}\label{prop:noexist-int}
$G_P$ cannot contain an interior vertex.
\end{proposition}

\begin{proof}
Assume not.
By Lemma \ref{lem:successive3}, any interior vertex of $G_P$ has degree at least $6$ in $\overline{G}_P$.
Hence $\overline{G}_P^+$ has either a good vertex of degree at most $4$ or a vertex of degree at most two by Proposition \ref{prop:int}.
In fact, there is a good vertex $x$ of degree $4$ by Lemma \ref{lem:4negative}.
Then $x$ is incident to $4$ successive positive edges of weight $4$.
(Examining the labels shows that there is no positive loop at $x$ as before.)
But this contradicts Lemma \ref{lem:successive3}.
\end{proof}

By Propositions \ref{prop:exist-int} and \ref{prop:noexist-int}, we have shown that $p\ge 3$ is impossible.

\subsection{Case 2: $p=2$}

There are $6$ possibilities (1)-(6) for $\overline{G}_P^+$ as stated in Lemma \ref{lem:2vertex}.
In particular, any vertex of $\overline{G}_P^+$ has degree at most $4$.

\begin{lemma}
$G_P$ has no interior vertex.
\end{lemma}

\begin{proof}
For an interior vertex, $4(t/2+2)\ge 5t$, giving $t\le 2$.
\end{proof}

\begin{lemma}\label{lem:negative-base}
Each vertex of $G_P$ is incident to a single negative loop in $\overline{G}_P$.
\end{lemma}

\begin{proof}
Assume that there is a vertex which is not incident to a negative loop.
If $G_P$ has no positive edges, then we have a contradiction by Lemma \ref{lem:key}.
Thus $G_P$ has a positive edge, then Proposition \ref{prop:key} implies $m=t=4$.
Hence $G_P$ has a good vertex of degree $4$.
But Lemma \ref{lem:degree4} works here.
Clearly, all negative loops at a vertex are parallel.
\end{proof}

Thus we can eliminate (1), (6) and Figure \ref{fig:2vertex}(2) of (3) in Lemma \ref{lem:2vertex}
as the possibility of $\overline{G}_P^+$ by Proposition \ref{prop:key}.
We now proceed to rule out the remaining possibilities of $\overline{G}_P^+$.

\begin{lemma}\label{lem:case2}
\textup{(2)} of Lemma \textup{\ref{lem:2vertex}} is impossible.
\end{lemma}

\begin{proof}
By Proposition \ref{prop:key}, $\widehat{T}$ is separating and $t=4$.
As in the proof of Lemma \ref{lem:degree2}, $\overline{G}_T$ has at least $8$ positive edges and at most $4$ negative edges.
By Lemma \ref{lem:parallel-max3}, $\overline{G}_T^+$ cannot have an isolated vertex or two vertices of degree one.
Thus $\overline{G}_T^+$ has two components, each of which is as shown in Figure \ref{fig:a-supp}(2).
Hence $\overline{G}_T$ has exactly $8$ positive edges, and so each vertex of $G_P$ is incident to $4$ negative loops (see the proof of Lemma \ref{lem:degree2}).
Let $k$ be the number of non-loop negative edges at $u_1$ in $G_P$.
Then $G_T$ contains $k$ positive $\{1,2\}$-edges.
Since a positive $\{1,2\}$-edge cannot be parallel to a positive level edge in $G_T$,
and since $G_T$ cannot contain an $S$-cycle, we see that $k=4$.
Thus each positive edge of $\overline{G}_P$ has weight $4$.

The core of the annulus support of $\overline{G}_P^+$ is separating or non-separating on $\widehat{P}$.
First, consider the case where it is separating.
Then the $4$ negative edges between $u_1$ and $u_2$ are divided into at most two families.
Since the associated permutation $\sigma$ to the negative loops at $u_1$ is the identity or $(13)(24)$,
the numbers of edges of those families are $\{4,0\}$ or $\{2,2\}$. 

For $\{4,0\}$, we may assume that $G_P$ has the labels as in Figure \ref{fig:cycle1}.
In this case, $\sigma$ is the identity.
Hence any vertex of $G_T$ is incident to a loop.
But there are only $4$ edges between $v_1$ and $v_2$.
This implies that any loop at $v_1$ is not level, a contradiction.

\nocolon\begin{figure}[htb]
\begin{center}
\includegraphics*[scale=0.32]{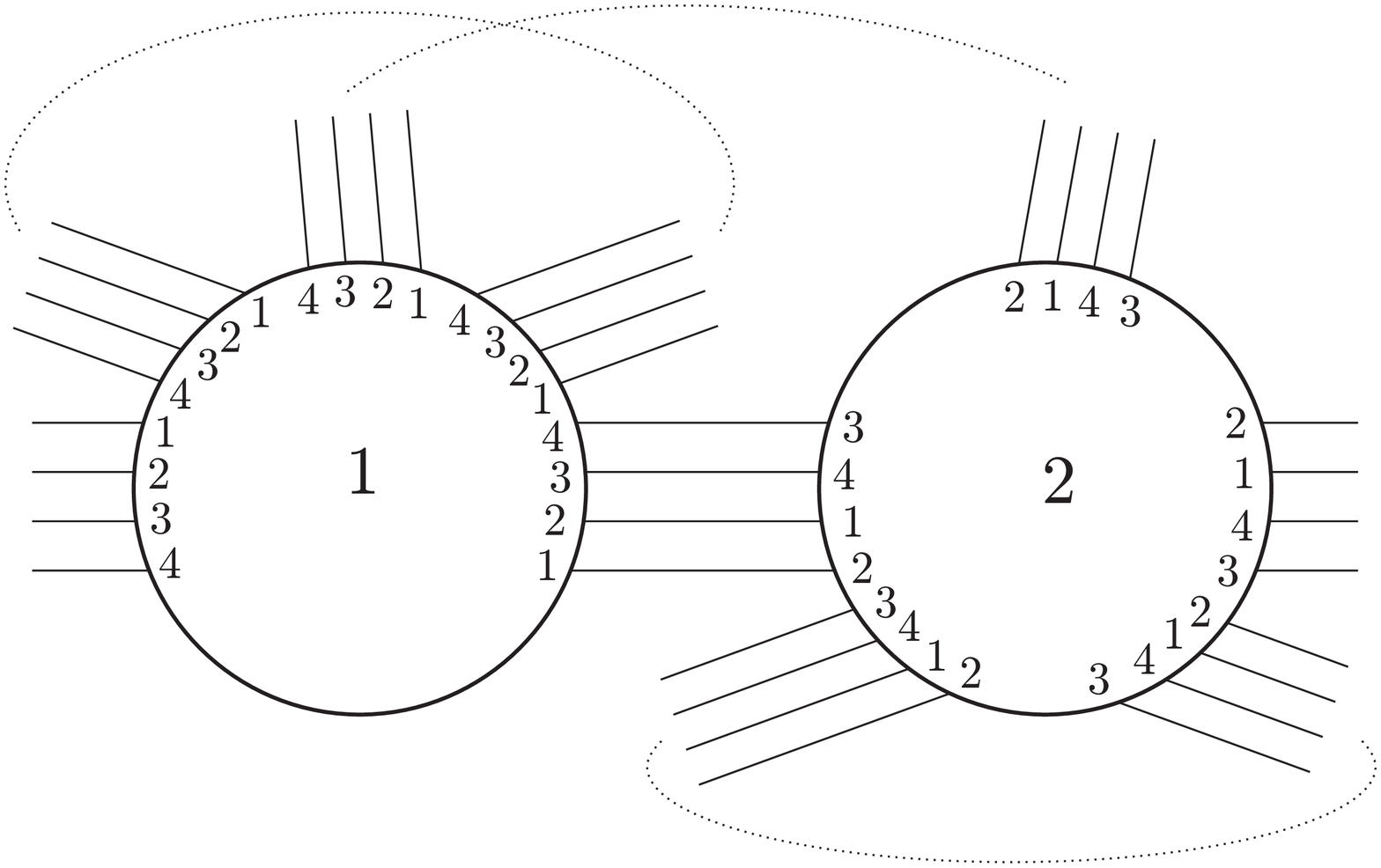}
\caption{}\label{fig:cycle1}
\end{center}
\end{figure}

For $\{2,2\}$, we may assume that $G_P$ has the labels as in Figure \ref{fig:cycle2}.
Then $\sigma=(13)(24)$.
Thus the negative loops at $u_1$ form two essential cycles on $\widehat{T}$.
Also, each vertex of $G_T$ is incident to a loop, and there are two parallel pairs between $v_i$ and $v_{i+1}$ for $i=1,3$.
Hence $G_T$ is uniquely determined.
But the arrangement of edges with label $1$ around $v_1$ contradicts Lemma \ref{lem:jumping}.
For, when we look at the two $1$-edges connecting $v_1$ with $v_3$, the jumping number must be two, but
when we look at the two $1$-edges connecting $v_1$ with $v_2$, it must be one.

\nocolon\begin{figure}[htb]
\begin{center}
\includegraphics*[scale=0.32]{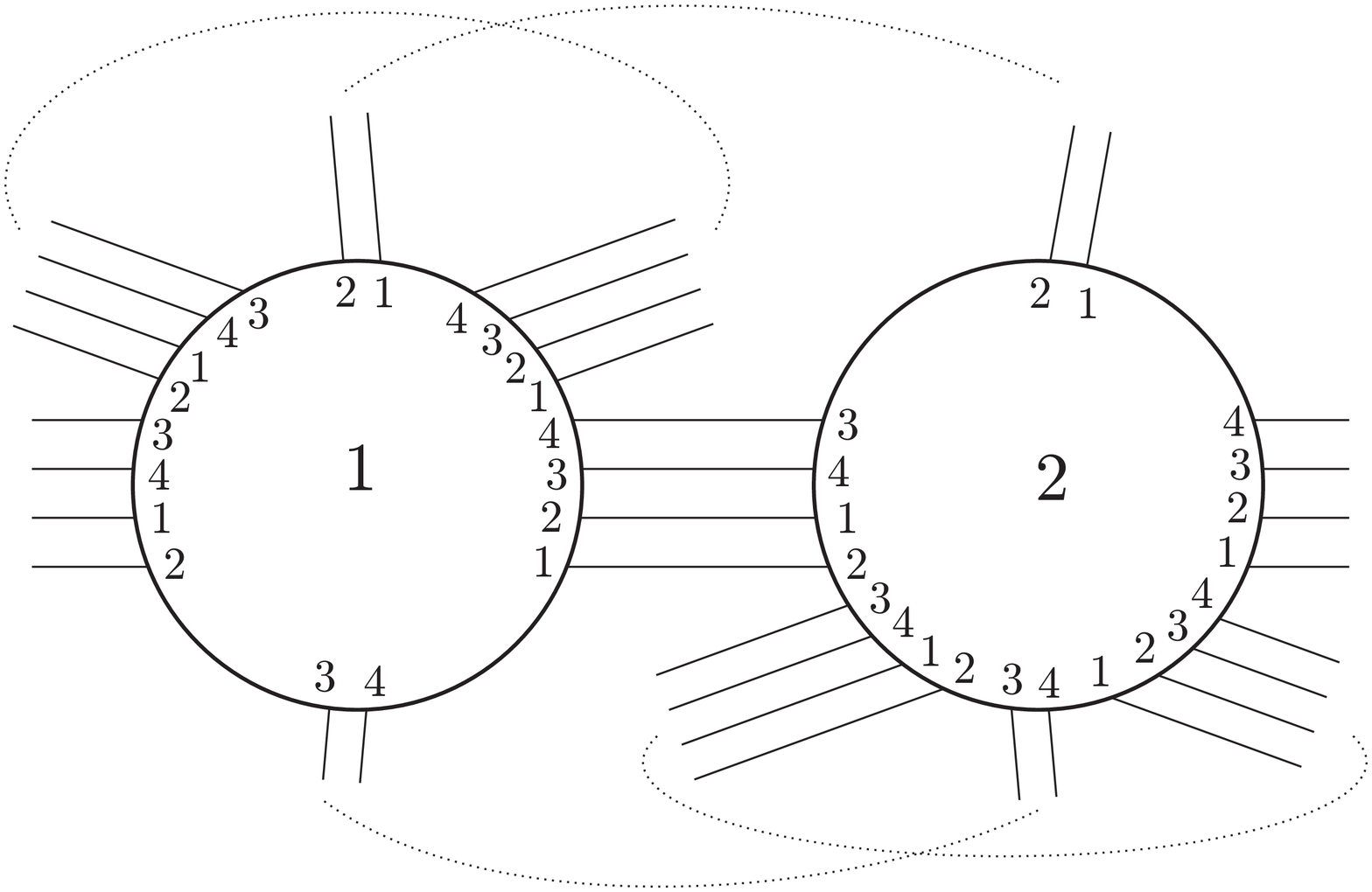}
\caption{}\label{fig:cycle2}
\end{center}
\end{figure}

Thus the core of the annulus support of $\overline{G}_P^+$ is non-separating on $\widehat{P}$.
We may assume that there are two $S$-cycles with label pair $\{1,2\}$ and $\{3,4\}$.
These $4$ edges form two essential cycles on $\widehat{T}$. 
Again, denote by $\sigma$ the associated permutation to the negative loops at $u_1$.
Then $\sigma$ is the identity or $(13)(24)$.
If $\sigma$ is the identity, then $G_T$ contains a graph as in Figure \ref{fig:cycle3}.
In $G_P$, there are $4$ negative edges between $u_1$ and $u_2$.
This means that $G_T$ has $4$ positive $\{1,2\}$-edges.
Since each vertex of $G_T$ has a level positive loop, these $4$ positive $\{1,2\}$-edges connect
$v_1$ and $v_3$ or $v_2$ and $v_4$, and furthermore three of them connect the same pair of vertices.
But then, there is a pair of parallel positive $\{1,2\}$-edges, which forms an $S$-cycle.

\nocolon\begin{figure}[htb]
\begin{center}
\includegraphics*[scale=0.54]{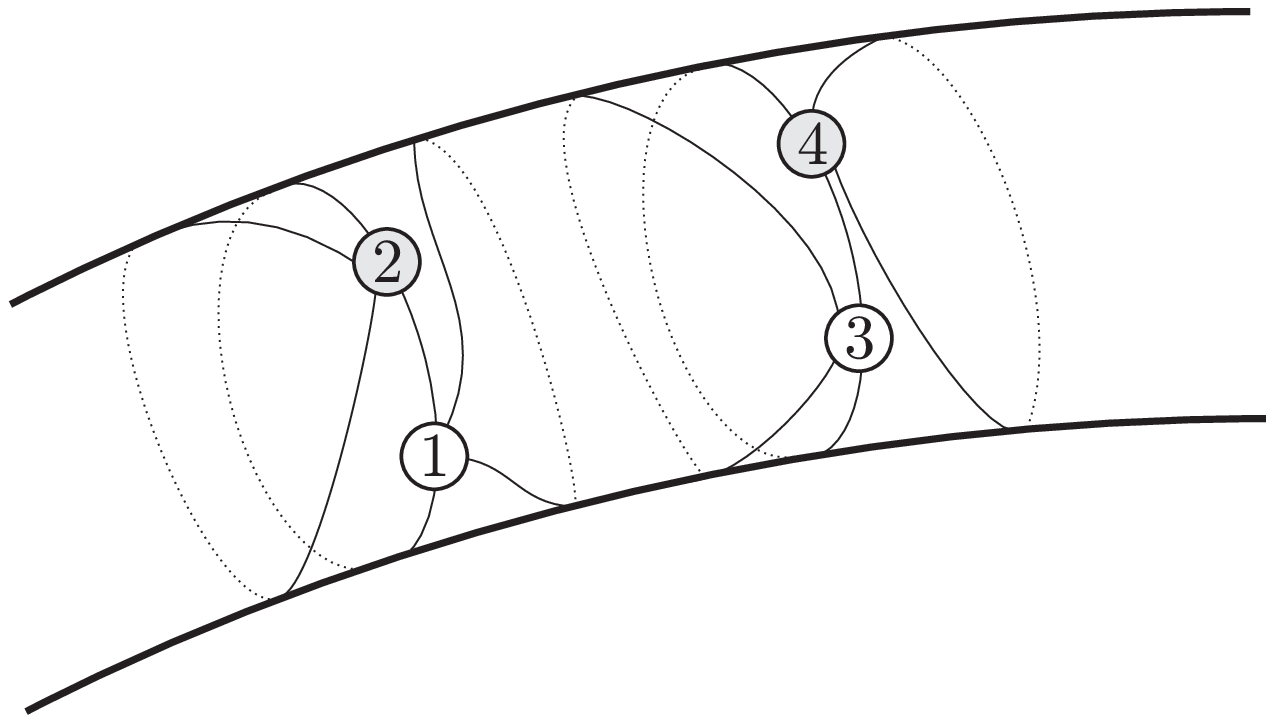}
\caption{}\label{fig:cycle3}
\end{center}
\end{figure}

\nocolon\begin{figure}[htb]
\begin{center}
\includegraphics*[scale=0.30]{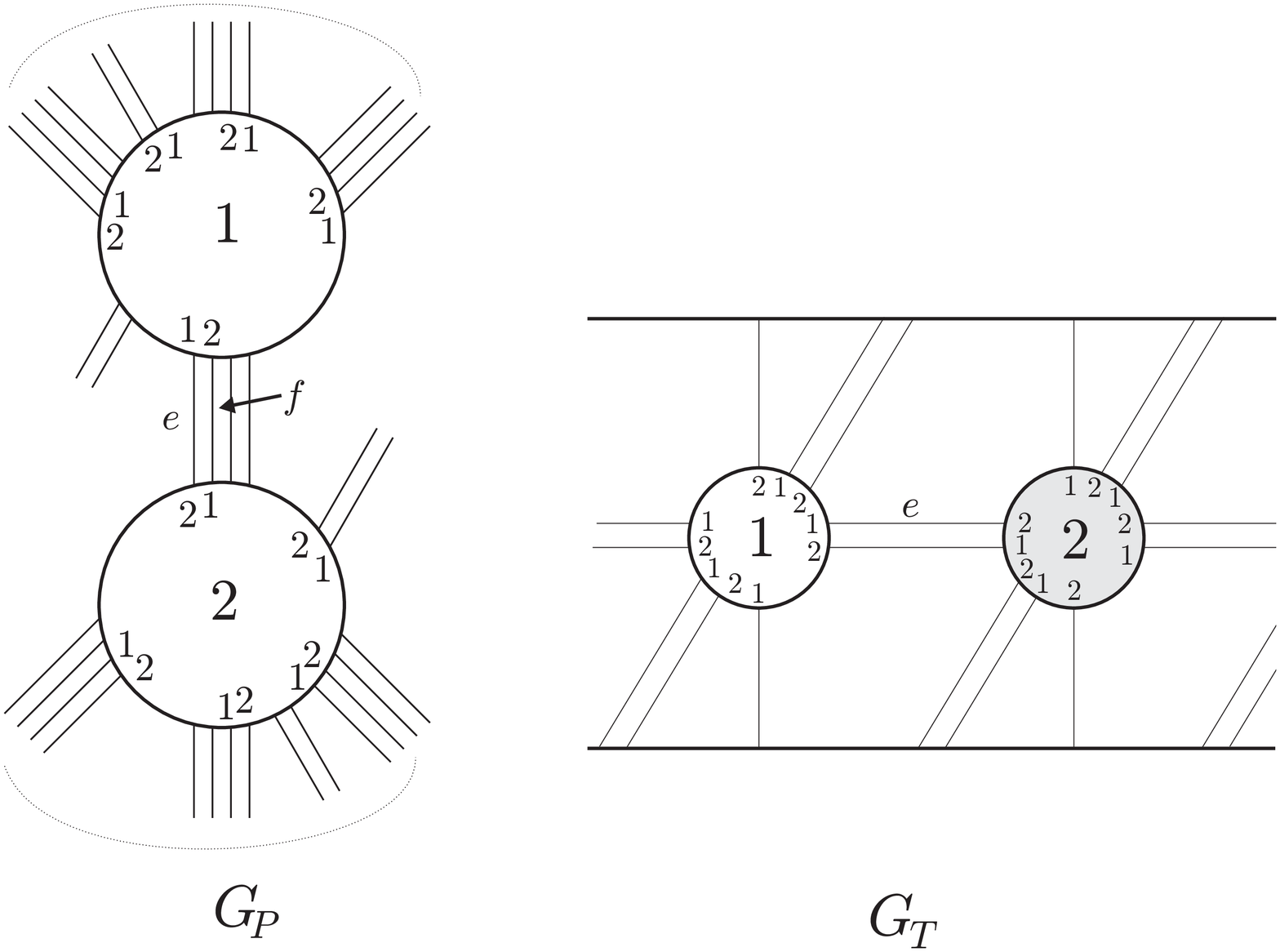}\vspace{-2mm}
\caption{}\label{fig:cycle4}
\end{center}
\end{figure}

Thus $\sigma=(13)(24)$, and then $G_P$ and hence $G_T$ are uniquely determined as in Figure \ref{fig:cycle4}.
At $u_1$, two occurrences of label $1$ at $\{1,2\}$ $S$-cycles are not consecutive among $5$ occurrences of label $1$.
But these points are consecutive at $v_1$.  Hence the jumping number is two.
Then the edge $e$ is located as in Figure \ref{fig:cycle4}.
Consider the location of the edge $f$ around $v_2$.
Then $f$ would be parallel to $e$, a contradiction.
\end{proof}

\begin{lemma}
\textup{(3)} of Lemma \textup{\ref{lem:2vertex}} is impossible.
\end{lemma}

\begin{proof}
Three cases are remaining.
For all cases, the core of the annulus support of $\overline{G}_P^+$ is separating on $\widehat{P}$ by
Lemma \ref{lem:negative-base}.
Also, when $\overline{G}_P^+$ is Figure \ref{fig:a-supp}(1) or (2), 
there is no non-loop negative edge in $G_P$.

Assume that $\overline{G}_P^+$ is as shown in Figure \ref{fig:a-supp}(1).
Hence $t=4$ by Lemma \ref{lem:key}(2), and then each vertex is incident to $4$ negative loops,
and each positive edge of $\overline{G}_P$ has weight $4$.
But there would be an extended $S$-cycle among positive loops, a contradiction.

Assume that $\overline{G}_P^+$ is as shown in Figure \ref{fig:a-supp}(2).
By Lemma \ref{lem:key}, $\widehat{T}$ is separating, and $t=4$ or $8$.
But, if $t=8$, then each positive edge of $\overline{G}_P$ has weight $6$ and each negative loop has weight $8$.
Then there would be an extended $S$-cycle among positive loops.
Thus $t=4$.

Let $n$ be the number of negative loops at $u_1$.
Then $n=3$ or $4$ by Lemmas \ref{lem:parallel-max-kb}(2) and \ref{lem:successive3}.
If $n=3$, then a negative loop at $u_1$ contradicts the parity rule.
Thus $n=4$.
By the same reason, $u_2$ is also incident to $4$ negative loops.
Then each vertex of $G_T$ is incident to two loops, which are level.

Let $l$ be the number of positive loops at $u_1$.
Then $l=2$, $3$ or $4$.

If $l=2$, we may assume that $G_P$ has the labels as in Figure \ref{fig:p2-2}(1).
Then $G_T$ has $4$ edges between $v_1$ and $v_2$.
But a loop at $v_1$ is level, a contradiction.

\nocolon\begin{figure}[htb]
\begin{center}
\includegraphics*[scale=0.32]{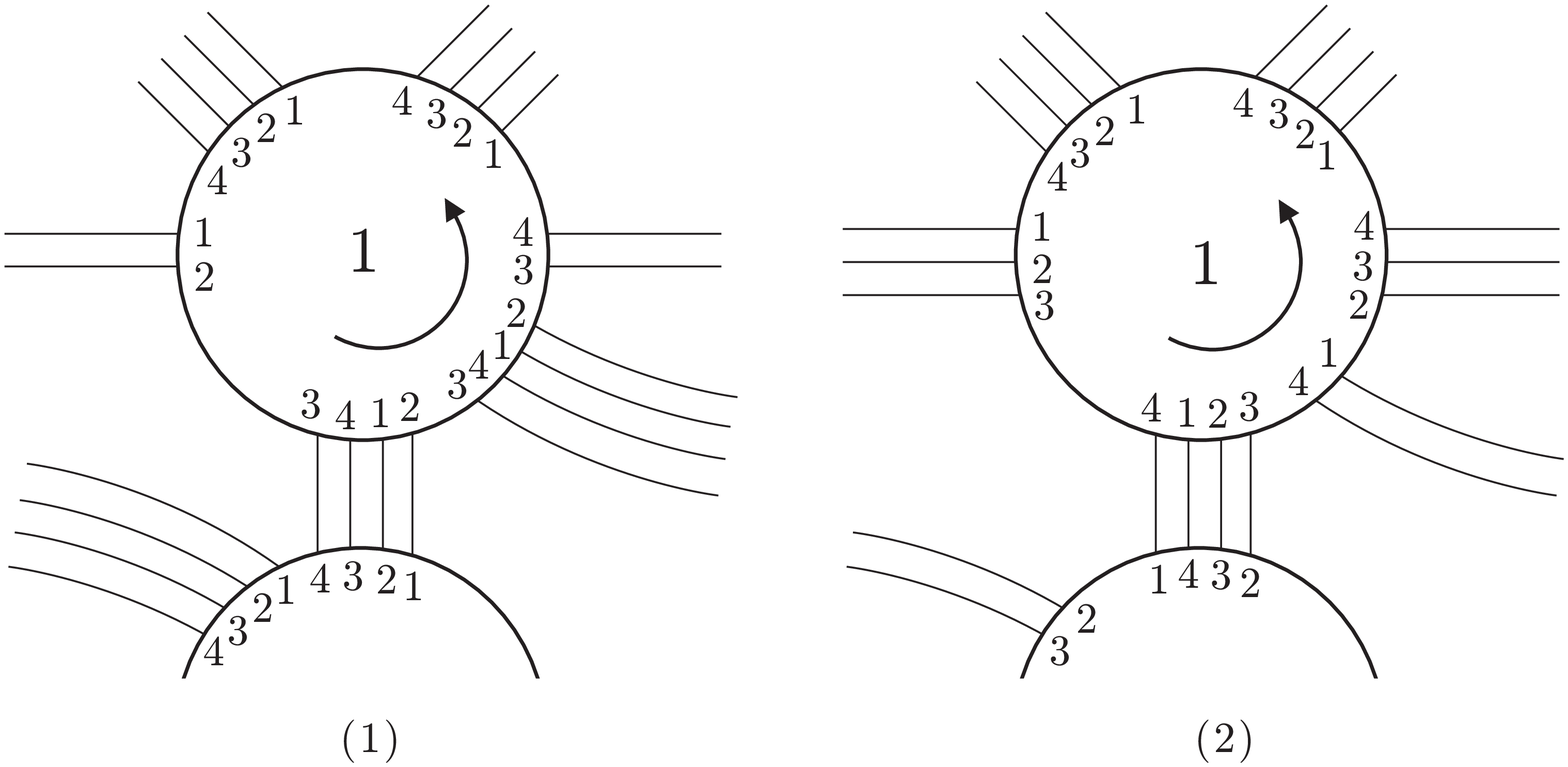}
\caption{}\label{fig:p2-2}
\end{center}
\end{figure}

If $l=3$, then two positive edges between $u_1$ and $u_2$ have weight $\{4,2\}$ or $\{3,3\}$. 
In the former, we may assume that $G_P$ has the labels as in Figure \ref{fig:p2-2}(2).
There are two $S$-cycles with label pair $\{2,3\}$.
In $G_T$, these $4$ edges are divided into two edge classes.
But such edge class cannot contain both a level edge and a non-level edge.
In the latter, we may assume that three positive loops at $u_1$ contain an $S$-cycle with label pair $\{2,3\}$ and a $\{1,4\}$-edge.
Although there are two possibilities for the labels at $u_2$, there is always a $\{2,3\}$-edge between $u_1$ and $u_2$.
Thus a similar argument to the former yields a contradiction.

If $l=4$, then there is an extended $S$-cycle among positive loops at $u_1$, a contradiction.

Finally, assume that $\overline{G}_P^+$ is as shown in Figure \ref{fig:2vertex}(1).
We may assume that $u_2$ is incident to a positive loop.
By Proposition \ref{prop:key} (applying to $u_1$), $\widehat{T}$ is separating and $t=4$.
Furthermore, $u_1$ is incident to $4$ negative loops, $4$ non-loop negative edges and
two families of $4$ parallel positive edges.
Thus we can assume that $G_P$ has the labels as in Figure \ref{fig:p2-3}.
Then the associated permutation to the negative loops is the identity, and that to the family $A$, say, of
non-loop negative edges is $(13)(24)$.
Let $e$ be the positive $\{2,3\}$-loop at $u_2$.

\nocolon\begin{figure}[htb]
\begin{center}
\includegraphics*[scale=0.32]{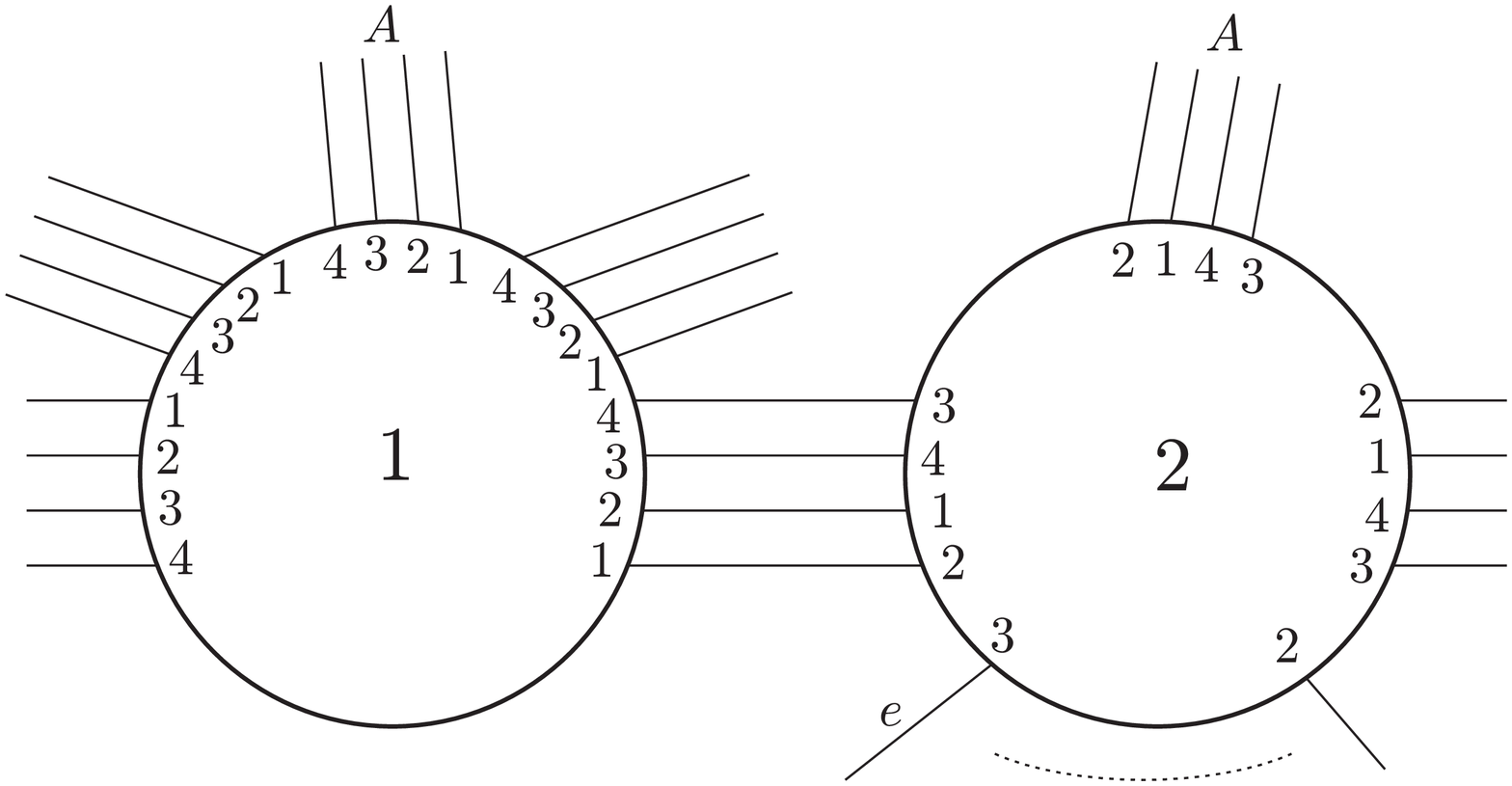}
\caption{}\label{fig:p2-3}
\end{center}
\end{figure}

The edges of $A$ form two essential cycles on $\widehat{T}$.
Put $e$ between $v_2$ and $v_3$.
Also, each vertex is incident to a loop.
Then we cannot locate the edges of an $S$-cycle with label pair $\{1,2\}$ so as to form an essential cycle.
\end{proof}

\begin{lemma}
If $\overline{G}_P^+$ is \textup{(4)} of Lemma \textup{\ref{lem:2vertex}}, then $\partial M$ consists of at most two tori.
\end{lemma}

\begin{proof}
Let $n$ and $k$ be the number of negative loops, non-loop negative edges, respectively, at $u_1$ in $G_P$.
Then $k+2n=5t$, and $k\le t$ by Lemma \ref{lem:key}.
Thus $n\ge 2t$.

Let $A$ be the family of negative loops at $u_1$, and let $\sigma$ be the associated permutation.
By Lemma \ref{lem:parallel-max-kb}, all vertices of $G_T$ have the same sign.
In fact, $\sigma$ has a single orbit \cite[Lemma 4.2]{Go2}.
Let $a_1,a_2,\dots,a_t,b_1$ be the successive $t+1$ edges in $A$, and let $D_1,D_2,\dots,D_t$ be
the disks between them.
Then the edges $a_1,a_2,\dots,a_t$ form an essential cycle on $\widehat{T}$, and
$a_2,a_3,\dots,a_t,b_1$ form a distinct essential cycle.
Let
\[ X=N(\widehat{T}\cup V_\beta\cup \bigcup_{i=1}^t D_i). \]
Then $\partial X$ is a torus $T'$, disjoint from $V_\beta$.
Thus either $T'$ is boundary parallel in $M(\beta)$, which implies that $\partial M$ is a union of two tori, or
$T'$ is compressible.
In the latter, $T'$ bounds a solid torus or is contained in a ball.
But if $T'$ lies in a ball, then $T$ would be compressible.
Hence $\partial M$ is a single torus.
\end{proof}

Remark that this case can be eliminated by a lengthy argument using a jumping number.

\begin{lemma}
\textup{(5)} of Lemma \textup{\ref{lem:2vertex}} is impossible.
\end{lemma}

\begin{proof}
Since each vertex of $G_P$ is incident to a negative loop by Lemma \ref{lem:negative-base},
the two loops of $\overline{G}_P^+$ are separating on $\widehat{P}$.
As in the proof of Lemma \ref{lem:case2}, $t=4$ and each vertex of $G_P$ is incident to $4$ negative loops and $4$ positive loops.
Then the family of $4$ loops contains an extended $S$-cycle, a contradiction.
\end{proof}

\subsection*{Acknowledgments}
We would like to thank the referee for his helpful comments.
Part of this work was done while we were visiting
The Mathematisches Forschungsinstitut Oberwolfach, Germany.
We are grateful for their hospitality and support.
This research was partially supported by Japan Society for the Promotion of Science,
Grant-in-Aid for Scientific Research (C), 16540071.


\Addresses\recd

\begin{thebibliography}{CGLS}


\bibitem{BGZ}
S. Boyer, C. McA. Gordon and X. Zhang,
\textit{Dehn fillings of large hyperbolic $3$-manifolds},
J. Differential Geom.\ \textbf{58} (2001), 263--308. 
  \MR{1913944}

\bibitem{BZ0}
S. Boyer and X. Zhang,
\textit{Reducing Dehn filling and toroidal Dehn filling},
Topology Appl.\ \textbf{68} (1996), 285--303.
  \MR{1377050}

\bibitem{CGLS}
M. Culler,  C. McA. Gordon, J. Luecke and P. Shalen,
\textit{Dehn surgery on knots},
Ann.\ of Math.\ \textbf {125} (1987), 237--300.
  \MR{0881270}

\bibitem{Go15}
C. McA. Gordon,
\textit{Dehn filling: a survey}, 
In \textit{Knot theory} (Banach Center Publ., 1998), 129--144.
  \MR{1634453}


\bibitem{Go2}
C. McA. Gordon,
\textit{Boundary slopes on punctured tori in $3$-manifolds},
Trans.\ Amer.\ Math.\ Soc.\ \textbf{350} (1998), 1713--1790. 
  \MR{1390037}

\bibitem{Go4}
C. McA. Gordon,
\textit{Small surfaces and Dehn fillings},
In \textit{Proceedings of the Kirbyfest}
(Berkeley, CA, 1998), 
\href{http://www.maths.warwick.ac.uk/gt/GTMon2/paper10.abs.html}
{Geom.\ Topol.\ Monogr. \textbf{2}, (1999),
177--199.} 
\MR{1734408}

\bibitem{GL2}
C. McA. Gordon and J. Luecke,
\textit{Dehn surgeries on knots creating essential tori, \rm{I}},
Comm.\ Anal.\ Geom.\ \textbf{3} (1995), 597--644.
  \MR{1371211}

\bibitem{GL3}
C. McA. Gordon and J. Luecke,
\textit{Toroidal and boundary-reducing Dehn fillings},
Topology Appl.\ \textbf{93} (1999), 77--90.
  \MR{1684214}

\bibitem{GL4}
C. McA. Gordon and J. Luecke,
\textit{Dehn surgeries on knots creating essential tori, \rm{II}},
Comm.\ Anal.\ Geom.\ \textbf{8} (2000), 671--725. 
  \MR{1792371}

\bibitem{GW}
C. McA. Gordon and Y. Q. Wu,
\textit{Toroidal and annular Dehn fillings}, 
Proc.\ London Math.\ Soc.\ \textbf{78} (1999), 662--700. 
  \MR{1674841}

\bibitem{HM}
C. Hayashi and K. Motegi,
\textit{Only single twists on unknots can produce composite knots}, 
Trans.\ Amer.\ Math.\ Soc.\ \textbf{349} (1997), 4465--4479. 
  \MR{1355073}

\bibitem{LOT}
S. Lee, S. Oh and M. Teragaito,
\textit{Reducing Dehn fillings and small surfaces}, 
preprint.

\bibitem{O}
S. Oh,
\textit{Reducible and toroidal $3$-manifolds obtained by Dehn fillings},
Topology Appl.\ \textbf{75} (1997), 93--104. 
  \MR{1425387}

\bibitem{O2}
S. Oh,
\textit{Reducing spheres and Klein bottles after Dehn fillings},
Canad.\ Math.\ Bull.\ \textbf{46} (2003), 265--267. 
  \MR{1981680}

\bibitem{P}
T. M. Price,
\textit{Homeomorphisms of quaternion space and projective planes in four space},
J. Austral.\ Math.\ Soc.\ \textbf{23} (1977), 112--128. 
  \MR{0436151}

\bibitem{T0}
M. Teragaito,
\textit{Creating Klein bottles by surgery on knots},
J. Knot Theory Ramifications \textbf{10} (2001), 781--794.
  \MR{1839702}

\bibitem{T}
M. Teragaito,
\textit{Distance between toroidal surgeries on hyperbolic knots in the $3$-sphere}, 
to appear in Trans.\ Amer.\ Math.\ Soc.


\bibitem{W}
Y. Q. Wu,
\textit{Dehn fillings producing reducible manifolds and toroidal manifolds},
Topology \textbf{37} (1998), 95--108.
  \MR{1480879}

\bibitem{W2}
Y. Q. Wu,
\textit{Sutured manifold hierarchies, essential laminations, and Dehn surgery},
J. Differential Geom.\ \textbf{48} (1998), 407--437.
  \MR{1638025}

\end{thebibliography}
\end{document}